%% file: main_with_bib.tex
\documentclass[a4paper,10pt]{article}

\usepackage[a4paper,left=3cm,right=3cm,top=2.5cm,bottom=2.5cm]{geometry}
\usepackage{hyperref}
\usepackage{subfig}
\usepackage{pgfplots}
\usepackage{pgfplotstable}
\usepackage{mathtools}
\usepackage{multicol}
\usepackage{comment}
\usepackage{booktabs}
\pgfplotsset{compat=1.5}
\usepackage{amssymb}
\usepackage{url}
\usepackage{bm}

\usepackage{pdfsync}
\usepackage{float}
\usepackage{tabularx}
\usepackage{enumerate}
\usepackage{array}
\usepackage{xspace}
\usepackage{siunitx}
\usepackage{authblk}

\usepackage{mathrsfs}
\usepackage{color, colortbl}
\usepackage{multirow}

\usepackage{amsthm}
\usepackage{amsmath}
\usepackage{stmaryrd}
\usepackage[ruled,vlined,linesnumbered]{algorithm2e}
\usepackage{graphicx}
\usepackage{booktabs}
\usepackage{cleveref}

\usepackage{tikz}
\usetikzlibrary{shapes,arrows,positioning,chains,fit,arrows.meta}

\usepackage[T1]{fontenc}

\DeclareMathOperator*{\argmax}{\arg\!\max}

\begin{document}

\title{Adaptive planning for risk-aware predictive digital twins}

\author[1]{Marco~Tezzele\footnote{marco.tezzele@austin.utexas.edu}}
\author[1]{Steven~Carr\footnote{stevencarr@utexas.edu}}
\author[1]{Ufuk~Topcu\footnote{utopcu@austin.utexas.edu}}
\author[1]{Karen~E.~Willcox\footnote{kwillcox@oden.utexas.edu}}

\affil[1]{Oden Institute for Computational Engineering and Sciences, University of Texas at Austin, Austin, 78712, TX, United States}

\maketitle

\begin{abstract}
This work proposes a mathematical framework to increase the robustness to rare events of digital twins modelled with graphical models. We incorporate probabilistic model-checking and linear programming into a dynamic Bayesian network to enable the construction of risk-averse digital twins. By modeling with a random variable the probability of the asset to transition from one state to another, we define a parametric Markov decision process. By solving this Markov decision process, we compute a policy that defines state-dependent optimal actions to take. To account for rare events connected to failures we leverage risk measures associated with the distribution of the random variables describing the transition probabilities. We refine the optimal policy at every time step resulting in a better trade off between operational costs and performances. We showcase the capabilities of the proposed framework with a structural digital twin of an unmanned aerial vehicle and its adaptive mission replanning.
\end{abstract}

\tableofcontents

\break

\section{Introduction}
\label{sec:intro}
When creating a digital twin (DT) it is important to ensure its robustness to different possible scenarios and rare event occurrences.
Moving from considering typical operational settings to more general scenarios that account for the involved risks will foster the development of robust and risk-averse digital twins. The aversion to risk is particularly important during the deployment of the DT, that is after the calibration phase is performed, to avoid interruption of the operations due to the occurrence of rare events.
To achieve reliability in an engineering system, it
  is important to account for the uncertainties characterizing the
  physical asset under consideration, in addition to the uncertainties
  in its operating environment. The main contribution of this work it
to present a method to increase the robustness of DTs by using a
probabilistic graphical model (PGM) formulation, parametric Markov
decision processes (MDPs), and risk measures to account for rare scenarios. The potential impact on real-world systems is twofold: reliability of DT deployments can be increased and operations can be optimized, resulting in lower operational costs and more effective predictive maintenance.

A DT is distinguished from traditional modeling and simulation through its personalization, that is, it is linked to the underlying physical asset or process, allowing it to represent the particular attributes of that specific asset or process. This DT representation is also dynamically evolving, continually incorporating data from the physical world, and refining the DT's predictive capabilities accordingly~\cite{aiaa2020digital, tuegel2011reengineering, grieves2017digital, rasheed2020digital, niederer2021scaling}.
Within the digital twin framework, we need to accomplish different
tasks such as data acquisition, solving inverse problems, building
computational models, accounting for uncertainty, and performing future
state predictions. In recent years, DTs have been proposed for diverse applications such as health monitoring and maintenance planning for spacecrafts~\cite{ye2020digital}, railway systems~\cite{arcieri2023bridging}, civil engineering structures~\cite{torzoni2024digital, morato2023inference}, and unmanned aerial vehicles~\cite{kapteyn2021probabilistic, kapteyn2022sensing, mcclellan2022physics}.

A common framework for modeling decision-making in the presence of uncertainty is based on Markov decision processes~\cite{puterman2014markov}. MDPs are widely used in reinforcement learning~\cite{sutton2018reinforcement}, robotics, and planning~\cite{russell2010artificial}. Oftentimes, the transitions in an MDP are unrealistically assumed to be precisely known. In this work, we leverage parametric MDPs with transition probabilities modeled as random variables with a given distribution. Parametric MDPs describe models where some probabilities are unknown but explicitly related. Such models form a natural pairing with a digital twin framework where information on one transition probability may inform future decisions. Probabilistic model-checking~\cite{kwiatkowska2017probabilistic} provides a general framework to verify whether an implemented policy satisfies a wide range of constraints.

The probabilistic graphical model we use is a dynamic Bayesian network with decision nodes, also called dynamic decision network~\cite{koller2009probabilistic, russell2010artificial, dean1989model,
  murphy2012machine}. The PGM encodes the relationships between the variables describing the physical asset and its digital counterpart. The edges in the PGM represent relationships between the connected variables through conditional probability distributions. The uncertainty is thus propagated in the entire graph, making it a powerful tool to account for risk. Among these conditional probability distributions the state transition probability distributions, which encode how the state evolves after taking a given action, are usually set a priori. In this paper, we propose a numerical method to dynamically update the state transition probabilities with a Bayesian approach, feeding the parametric MDP with the new information, and making the graphical model more interconnected with the physical asset. Better estimates translate into better policies and a more effective digital twin. 

\begin{figure}[ht]
  \centering
  \includegraphics[width=1.\textwidth, trim=170 40 90 40, clip]{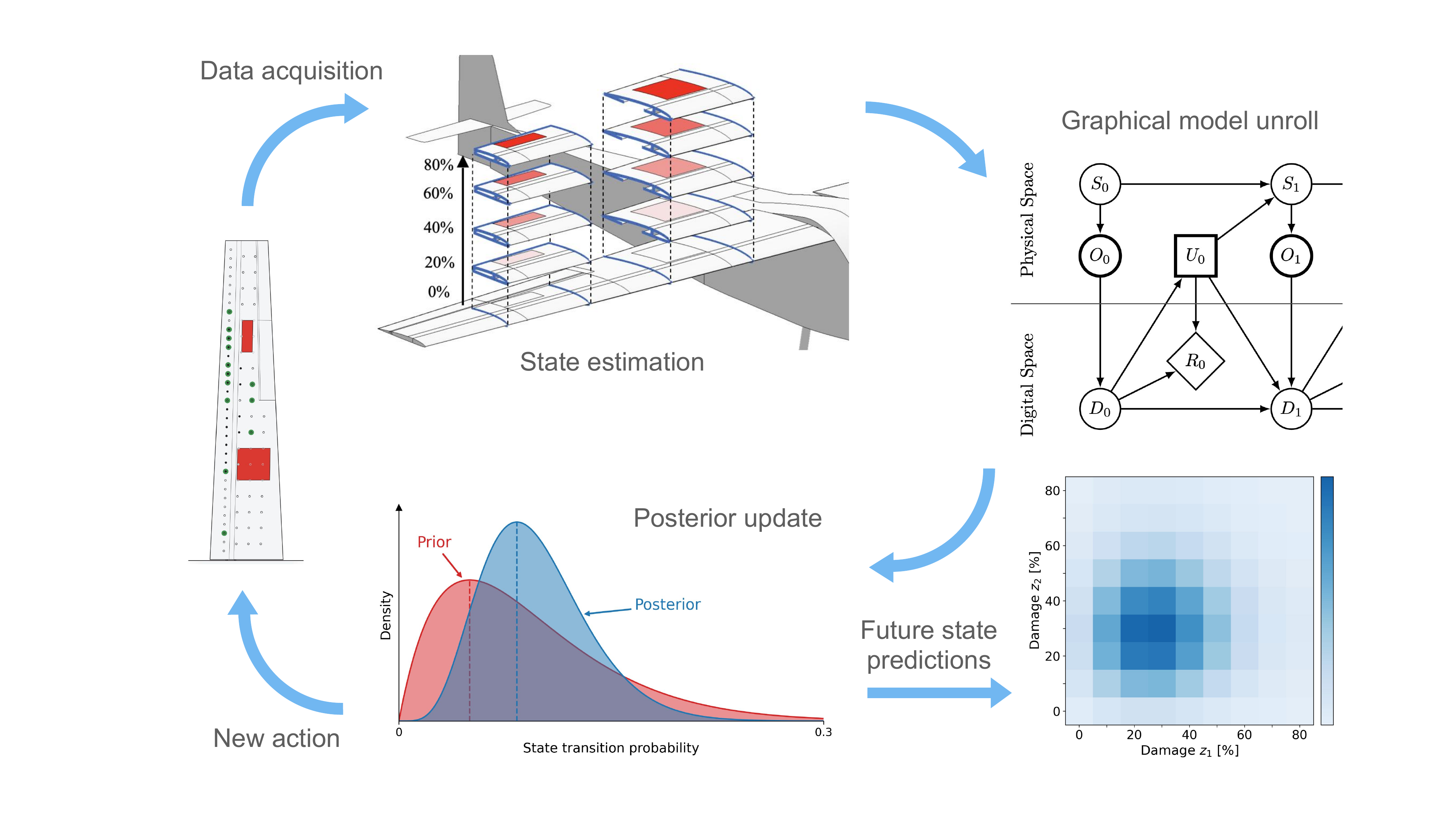}
  \caption{Abstract representation of the information flow within the DT formulation, from sensor data acquisition to the MDP's optimal policy update. Starting from strain sensors placed on the wing, we estimate the digital state connected with the structural integrity of the UAV. We then unroll the graphical model with all the variables for a new time step and we update the posterior estimates of the transition probability. We recompute the new optimal policy and we issue a new action to perform. This policy is also used to compute future actions in the prediction phase happening in the digital space.}
  \label{fig:pipeline_scheme}
\end{figure}

A sketch of the end-to-end pipeline is depicted in Figure~\ref{fig:pipeline_scheme}, where we consider a specific physical asset for the sake of clarity. The computational flow starts with the data acquisition from sensors placed on the unmanned aerial vehicle (UAV) wing. The noisy data are used to estimate the digital state and the PGM is populated with another time step. The unrolling is realized by duplicating the nodes and edges of the graph, at a given time, for the next time step. The policy is then updated by solving the MDP using the new state transition probability, and the optimal action is issued, continuing the computational cycle. Future state predictions can be computed for any time step, using the current estimate of the transition probability.

The paper is organized as follows. In Section~\ref{sec:pgm} we describe the PGM and its predictive capabilities. Section~\ref{sec:policy} presents the underlying parametric MDP and how the optimal policy is updated using different risk measures. The numerical results with a test mission and two different scenarios are explained in Section~\ref{sec:results}. Finally, we draw conclusions in Section~\ref{sec:the_end}.

\section{Probabilistic graphical model formulation for digital twins}
\label{sec:pgm}
This section introduces the PGM and how the variables are connected with the MDP. We follow the approach of Ref.~\cite{kapteyn2021probabilistic} for the PGM formulation.

\begin{figure}[ht]
\centering
\input{figures/digital_twin_graph.tex}
\caption{Representation of the dynamic decision network used to encode
  the relationships between the physical and the digital spaces for the first 3 time steps. Square nodes denote actions, while diamond nodes denote the reward function. Bold outlines represent observed (deterministic) quantities, while thin outlines stand for estimated (random) variables.}\label{fig:digital_twin_graph}
\end{figure}
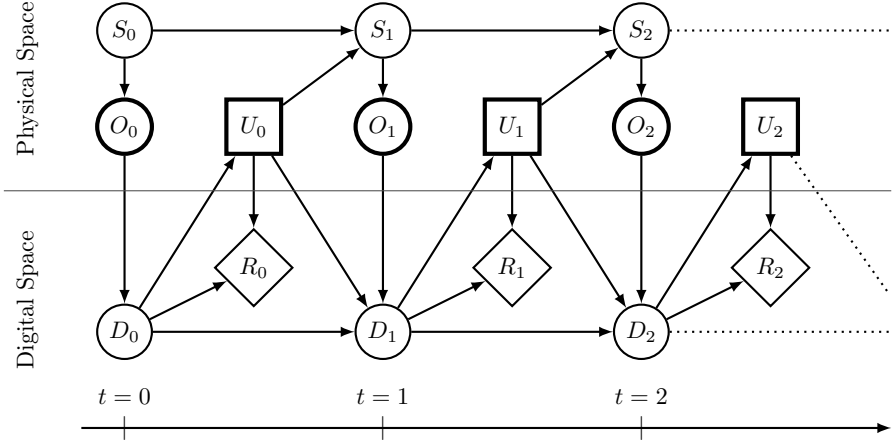

Figure~\ref{fig:digital_twin_graph} depicts the encoding of the
interactions between the physical and the digital assets through a
dynamic decision network, starting from $t = 0$ to the current time $t = t_c =
2$. We indicate with $t$ a generic non-dimensional discretized time step.  Circle nodes in the graph represent random variables, diamonds indicate output functions, while square nodes denote the decisions taken to enact specific
actions. Following the edges we track the
dependencies between variables. Notice how we can unroll the directed
graphical model for every time step just by repeating its core
structure. We consider the following time-dependent random variables: the unknown underlying physical state $S$, the observed sensor measurements $O$, the digital state $D$, the actions $U$ (also called controls) that can be issued, and the reward $R$. The direct solid edges encode the temporal dependencies between the variables through conditional probability distributions. 

Decision-making is achieved by leveraging MDPs. Let $\mathcal{S}$ be the set of possible states, $\mathcal{A}$ be the set of possible actions, $\mathcal{P}$ be a set of transition probabilities, and $\mathcal{R}$ be the reward set. In the remainder of the work we assume all these sets are finite. A finite Markov decision process is defined by the tuple $(\mathcal{S}, \mathcal{A}, \mathcal{P}, \mathcal{R})$. The state transition probability $p: \mathcal{S} \times \mathcal{S} \times
\mathcal{A} \to [0, 1]$ describes the probability of transitioning from
the state $s$ to $s^\prime$ at time $t-1$ after taking the action
$u$. It is defined as:
\begin{equation}
p(s^\prime \mid s, u) := \text{Pr} \{ S_t = s^\prime \mid S_{t-1} = s, U_{t-1} = u \}.
\end{equation}
Analogously the reward function $R$ is defined as:
\begin{equation}
R(s^\prime, s, u) := \text{Pr} \{ r_t \mid S_t = s^\prime, S_{t-1} = s,
U_{t-1} = u \},
\end{equation}
where $r_t$ is the reward at time $t$. The probability of taking an action $u$ at time $t$ is given by a policy $\pi_t: \mathcal{S} \times
\mathcal{A} \to [0, 1]$ defined as:
\begin{equation}
\pi_t(s, u) := \text{Pr} \{ U_t = u \mid S_{t-1} = s \}.
\end{equation}
To issue an actual action $u_t$ we take the $\argmax_{u \in \mathcal{A}} \pi_t(s, u)$. A classical MDP satisfies the Markov property~\cite{puterman2014markov}, that is 
\begin{equation}
\begin{split}
\text{Pr} \{ S_t = s^\prime \mid S_{t-1} = s_{t-1}, U_{t-1} = u_{t-1},
\dots, S_1 = s_1, U_1 = u_1 \} = \\ \text{Pr} \{ S_t = s^\prime \mid S_{t-1} = s_{t-1}, U_{t-1} = u_{t-1} \}.
\end{split}
\end{equation}
The goal is to leverage the resolution of an MDP to find the optimal sequence of actions to issue to maximize the expected reward. The MDP is thus a key component to provide optimal decision-making capabilities to the PGM. For more general cases where the Markov property does not hold, the topology of the graph should reflect that by adding more edges connecting all the necessary time instants.

We now show how uncertainty is propagated through the graph from data acquisition to future state predictions. We denote by $p(d_t)$ the probability distribution of $D_t=d_t$, for any possible digital state $d_t$, with similar notation employed for the other random variables comprising the graphical model. We have $S_t \sim p(s_t)$ for the physical state, $O_t \sim p(o_t)$ for the observational data, $U_t \sim p(u_t)$ for the action, and $R_t \sim p(r_t)$ for the reward.
We can factorize the joint distributions of some of the variables conditioned on the observed ones, as in the following:
\begin{equation*}
p(D_0, \dots, D_{t_c}, R_{0}, \dots, R_{t_c} \mid o_0, \dots, o_{t_c}, u_0, \dots, u_{t_c})
\propto \prod_{t=0}^{t_c}\Bigl[\phi_t^{\text{data}} \phi_t^{\text{dynamics}} \phi_t^{\text{reward}} \Bigr],
\label{eq:factorization_graph}
\end{equation*}
where the factors are:
\begin{align}
\phi_t^{\text{data}} &= p(O_t = o_t \mid D_t), \\
\phi_t^{\text{dynamics}} &= p(D_t \mid D_{t-1}, U_{t-1} = u_{t-1}), \\
\phi_t^{\text{reward}} &= p(R_t \mid D_t, U_t = u_t).
\end{align}
The factor $\phi_t^{\text{data}}$ encodes the assimilation of the sensor
  measurements. With $\phi_t^{\text{dynamics}}$ we factorize the
belief about the digital state at the current time step, given the
digital state at the previous time step and the last enacted control. Finally $\phi_t^{\text{reward}}$ encodes the evaluation of the reward function which is a performance measure. The state space and the action space are discrete so the conditional probabilities involving these spaces are encoded in the graph through conditional probability tables (CPTs). These CPTs are usually set a priori, extrapolated from historical mission logs or manufacturer's blueprints of the physical asset. The structure of the CPTs describing $\phi_t^{\text{dynamics}}$, also called transition probability matrices, is further analyzed in Section~\ref{sec:mdp}.

\begin{figure}[h!]
\centering
\input{figures/prediction_graph.tex}
\caption{Dynamic Bayesian network used for the prediction of the digital state evolution with the associated uncertainty. We assume to have the estimation of $D_{t_c}$ given the observations from the physical asset $O_{t_c}$, the previous issued action $U_{t_c-1}$, and the previous digital state $D_{t_c-1}$. From the current time $t_c$ we show the graph used to predict the next 3 actions and digital states, using the policy $\pi_{t_c}$.}\label{fig:prediction_graph}
\end{figure}
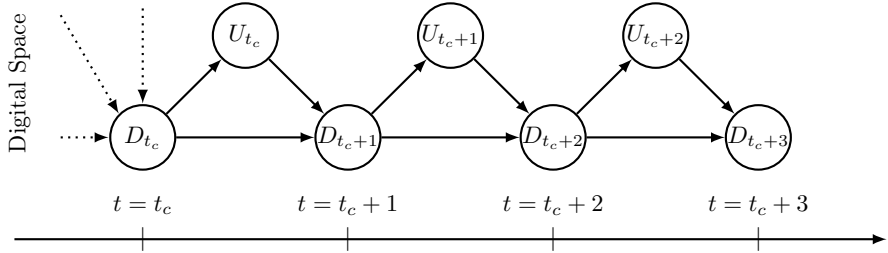

The prediction for future time steps is performed using just a subset
of the variables comprising the DT framework, as depicted in
Figure~\ref{fig:prediction_graph}. To decide which action to perform we query the policy $\pi$, mapping the current digital state to a distribution of actions. The action with the highest density is the one that maximizes the total
future rewards. After updating the transition probability estimates in the graph at every time step, we need to recompute the corresponding optimal policy. The policy used for all the times $t \geq t_c$ is denoted with $\pi_{t_c}$ and incorporates the most updated state transition probability estimates for all the possible actions. Notice how in Figure~\ref{fig:prediction_graph} the actions $U_t$ are depicted as random variables and not anymore as decision nodes since we propagate the uncertainty throughout all the prediction stages. We have indeed that $U_{t_c+i} \sim \pi_{t_c}(D_{t_c+i})$ for $i=0, 1, 2$.

The target belief state is extended as in the following:
\begin{align}
&p(D_0, \ldots, D_{t_p}, R_{0}, \dots, R_{t_c}, U_{t_c + 1}, \ldots, U_{t_p} \mid o_0, \ldots, o_{t_c}, u_0, \ldots, u_{t_c}) \nonumber \\
&\propto 
 \prod_{t=0}^{t_c} 
 \Bigl[ \phi_t^{\text{data}} \phi_t^{\text{dynamics}} \phi_t^{\text{reward}} \Bigr] \prod_{t=t_c + 1}^{t_p} \Bigl[ \phi_t^{\text{step}} \phi_t^{\text{control}} \Bigr],
\label{eq:factorization_pred}
\end{align}
where $\phi_t^{\text{control}} = p(U_t \mid D_t)$ is given by the last computed policy $\pi_{t_c}$, and $\phi_t^{\text{step}} = p(D_t \mid D_{t-1}, U_{t-1})$.

\section{Online adaptivity of risk-aware parametrized policies}
\label{sec:policy}
Section~\ref{sec:mdp} introduces parametric Markov decision processes defined by random variables and how the state transition probability is refined after each data acquisition step. Section~\ref{sec:risk} presents how risk measures are incorporated into the PGM formulation.

\subsection{Parametric Markov decision processes}
\label{sec:mdp}

For every possible action, there is a corresponding transition probability matrix encoding the conditional probabilities for all the combinations of states. The transition probability matrix $P$ associated to the action $u \in \mathcal{A}$  is defined as
\begin{equation}
P(u) := (p_{ij})_{i, j = 1}^n = (p(s_j \mid s_i, u))_{i, j = 1}^n, \qquad \text{for } \, s_i, s_j \in \mathcal{S},
\end{equation}
that is:
\begin{equation}\label{eq:transition_matrix}
P(u) =
  \begin{pmatrix}
    p(s_1 \mid s_1, u) & p(s_2 \mid s_1, u) & \dots & p(s_n \mid s_1, u)  \\
    p(s_1 \mid s_2, u) & p(s_2 \mid s_2, u) & \dots & p(s_n \mid s_2, u)  \\
    \vdots & \vdots & \ddots & \vdots \\
    p(s_1 \mid s_n, u) & p(s_2 \mid s_n, u) & \dots & p(s_n \mid s_n, u)  \\
  \end{pmatrix} ,
\end{equation}
where $p_{ij} \geq 0$, and $\sum_{j=1}^n p_{ij} = 1$, for all $i \in
[1, \dots, n]$. This matrix can take different forms depending on
the underlying MDP.

We assume there exists a set of actions,
denoted by $\mathcal{A}^D
\subset \mathcal{A}$, which produces a deterministic outcome such as moving towards a known direction, for example. We define
$\mathcal{A}^D$ as
\begin{equation}
\mathcal{A}^D := \{ u \in \mathcal{A} \mid \forall s \in \mathcal{S}, \,\, \exists \,
s^\prime \in \mathcal{S} \,\text{ s.t. }\, p(s^\prime \mid s, u) = 1 \}.
\end{equation}
We denote its complementary set with $\mathcal{A}^N = \mathcal{A}
\setminus \mathcal{A}^D$ and we define it as
\begin{align}
  \label{eq:nondet_actions}
\mathcal{A}^N := \{ &u \in \mathcal{A} \mid \forall s \in \mathcal{S} \setminus \{ \bar{s} \}, \,\, \exists \,
s^\prime, s^{\prime\prime} \in  \mathcal{S}, \,\, s^\prime \neq
s^{\prime\prime}, \,\, \exists \, q \in (0, 1) \,\text{ s.t. }\,\\
&p(s^\prime \mid s, u) = q, \,\, p(s^{\prime\prime} \mid s, u) = 1-q \},\nonumber
\end{align}
where $\bar{s}$ is the only state for which $p(\bar{s} \mid \bar{s}, u) = 1$.

We restrict to the case for which only two outcomes
are possible, i.e., the state varies smoothly without abrupt changes that lead to the transition to many possible states. This assumption is not too restrictive because we can always decrease the acquisition time step up to a certain computational threshold. Moreover, $s^\prime$ and $s^{\prime\prime}$ in Eq.~\eqref{eq:nondet_actions} will always be two consecutive discretized states.

\begin{figure}[ht]
\centering
\begin{tikzpicture}[every node/.style={minimum size=.5cm-\pgflinewidth, outer sep=0pt}]
    \node at (1.5, 3.3) {$P(u) \,\,\, \forall u \in \mathcal{A}^D$};
    \node at (-0.3, 1.5) {$s$};
    \node at (1.5, -0.3) {$s^\prime$};
    \draw[step=0.5cm,color=black] (0,0) grid (3,3);
    \node[draw=darkgray, line width=0.12mm, fill=green!70!cyan!40] at (2.25, 0.25) {};
    \node[draw=darkgray, line width=0.12mm, fill=green!70!cyan!40] at (0.25, 0.75) {};
    \node[draw=darkgray, line width=0.12mm, fill=green!70!cyan!40] at (1.75, 1.25) {};
    \node[draw=darkgray, line width=0.12mm, fill=green!70!cyan!40] at (0.75, 1.75) {};
    \node[draw=darkgray, line width=0.12mm, fill=green!70!cyan!40] at (2.75, 2.25) {};
    \node[draw=darkgray, line width=0.12mm, fill=green!70!cyan!40] at (1.25, 2.75) {};
\end{tikzpicture} \qquad
\begin{tikzpicture}[every node/.style={minimum size=.5cm-\pgflinewidth, outer sep=0pt}]
    \node at (1.5, 3.3) {$P(u; q), \,\,\, u \in \mathcal{A}^N$};
    \node at (-0.3, 1.5) {$s$};
    \node at (1.5, -0.3) {$s^\prime$};
    \draw[step=0.5cm,color=black] (0,0) grid (3,3);
    \node[draw=darkgray, line width=0.12mm, fill=green!70!cyan!40] at (2.75, 0.25) {};
    \node[draw=darkgray, line width=0.12mm, fill=red!60] at (2.25, 0.75) {};
    \node[draw=darkgray, line width=0.12mm, fill=red!60] at (1.75, 1.25) {};
    \node[draw=darkgray, line width=0.12mm, fill=red!60] at (1.25, 1.75) {};
    \node[draw=darkgray, line width=0.12mm, fill=red!60] at (0.75, 2.25) {};
    \node[draw=darkgray, line width=0.12mm, fill=red!60] at (0.25, 2.75) {};
    \node[draw=darkgray, line width=0.12mm, fill=cyan!60] at (2.75, 0.75) {};
    \node[draw=darkgray, line width=0.12mm, fill=cyan!60] at (2.25, 1.25) {};
    \node[draw=darkgray, line width=0.12mm, fill=cyan!60] at (1.75, 1.75) {};
    \node[draw=darkgray, line width=0.12mm, fill=cyan!60] at (1.25, 2.25) {};
    \node[draw=darkgray, line width=0.12mm, fill=cyan!60] at (0.75, 2.75) {};

    \draw [dashed] (4, -0.5) -- (4, 3.7);
  \end{tikzpicture} \qquad
  \begin{tikzpicture}[every node/.style={minimum size=.5cm-\pgflinewidth, outer sep=0pt}]
    \node at (1.5, 3.3) {$P(u; q(s)), \,\,\, u \in \mathcal{A}^N$};
    \node at (-0.3, 1.5) {$s$};
    \node at (1.5, -0.3) {$s^\prime$};
    \draw[step=0.5cm,color=black] (0,0) grid (3,3);
    \node[draw=darkgray, line width=0.12mm, fill=green!70!cyan!40] at (2.75, 0.25) {};
    \node[draw=darkgray, line width=0.12mm, fill=yellow!80] at (2.25, 0.75) {};
    \node[draw=darkgray, line width=0.12mm, fill=purple!40] at (1.75, 1.25) {};
    \node[draw=darkgray, line width=0.12mm, fill=red!60] at (1.25, 1.75) {};
    \node[draw=darkgray, line width=0.12mm, fill=blue!30] at (0.75, 2.25) {};
    \node[draw=darkgray, line width=0.12mm, fill=magenta!40] at (0.25, 2.75) {};
    \node[draw=darkgray, line width=0.12mm, fill=blue!20!cyan!30] at (2.75, 0.75) {};
    \node[draw=darkgray, line width=0.12mm, fill=orange!80] at (2.25, 1.25) {};
    \node[draw=darkgray, line width=0.12mm, fill=black!30] at (1.75, 1.75) {};
    \node[draw=darkgray, line width=0.12mm, fill=green!40] at (1.25, 2.25) {};
    \node[draw=darkgray, line width=0.12mm, fill=cyan] at (0.75, 2.75) {};
\end{tikzpicture}
\caption{Illustration of the types of transition matrices considered in this work. On
  the left panel the matrix for the actions in $\mathcal{A}^D$, where
  green stands for $1$ and white for $0$. In the central panel we have
  the upper bidiagonal with fixed transitions for actions in
  $\mathcal{A}^N$, where red stands for $1-q$ and blue for $q$. On the
  right panel the generic upper bidiagonal for
  transition probabilities varying between states.}
    \label{fig:matrix_structures}
\end{figure}
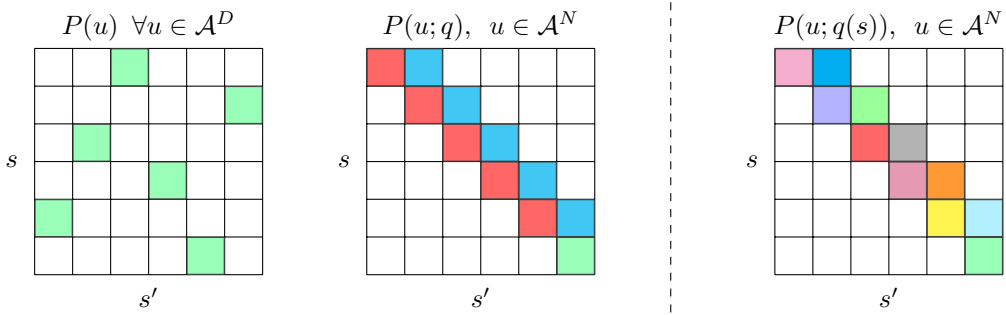

Let us consider an action $u^D \in \mathcal{A}^D$,
then we have that $P(u^D)$ is a permutation of the identity matrix, which
is known a priori. If, instead, we consider an action $u^N \in
\mathcal{A}^N$, then we can parametrize $P(u^N)$ with the parameter $q$ introduced in Eq.~\eqref{eq:nondet_actions}. We introduce $P(u^N; q)$, an upper bidiagonal matrix
\footnote{In general the number of diagonals is equal to $2^d$, where $d$ is the dimension of the state $s$.}, which has the following form:
\begin{equation}\label{eq:parametrized_transition_matrix}
P(u^N; q) =
  \begin{pmatrix}
    1-q & q & 0 & \dots & 0 & 0  \\
    0 & 1-q & q & \dots & 0 & 0  \\
    \vdots & \vdots & \vdots & \ddots & \vdots & \vdots \\
    0 & 0 & 0 & \dots & 1-q & q  \\
    0 & 0 & 0 & \dots & 0 & 1  \\
  \end{pmatrix} .
\end{equation}
We refer to this case as state-invariant transition, because the probability of transitioning from the current state to the contiguous discretized state does not depend on the current state. In particular, if $s$ is a scalar, the superdiagonal elements are equal to $q$, while the main diagonal elements are $1-q$, except the last one which is $1$.  Figure~\ref{fig:matrix_structures} illustrates the different transition probability matrices arising from our setting. The general
case with a state-dependent transition, denoted as $P(u;
q(s))$ in Figure~\ref{fig:matrix_structures}, will be considered in
future studies.

In general, $q$ is unknown and possibly varying in time. When the parameter $q$ is incorporated to the transition probabilities, we can model the system as a parametric MDP, which describes sets of sets on the transition probabilities. We model $q$ with a random variable $Q$. By accounting for the unknown nature of $q$, we may compute robust policies covering the worst case transitions~\cite{CubuktepeJJKT22}. As we gather more information on the true value of $q$ the policy will dynamically improve. We leverage the PGM formulation and the acquired observations to refine our belief on $q$ at every time steps. Other approaches to describe uncertainty in parametric MDPs use intervals~\cite{delahaye2011decision, givan2000bounded, puggelli2013polynomial}, robust MDPs~\cite{wiesemann2013robust} which consider nonconvex sets for the transition probabilities~\cite{nilim2005robust}, or temporal specifications~\cite{wolff2012robust}.

Suppose that the random variable $Q$ describing the state transition
probability of our parametric MDP is beta-distributed with left parameter
$\alpha>1$ and right parameter $\beta>1$, that is $Q \sim \mathcal{B}e(\alpha,
\beta)$. Let us also suppose that $X = (X_1, X_2, \dots)$ is a
sequence of indicator random variables such that given $Q = q \in (0,
1)$, $X$ is a conditionally independent sequence with
\begin{equation}
\text{Pr} \{ X_i = 1 \mid Q = q \} = q, \qquad \forall i \in \mathbb{N}^+,
\end{equation}
then $X$ is called a beta-Bernoulli process with parameters $\alpha$
and $\beta$. $X_i$ represents the outcome of the $i$-th trial,
with $1$ denoting a success and $0$ a failure. In our particular case
$X_i = 1$ if the state at time $i$ is different from the state at time $i-1$, and $X_i = 0$ otherwise.

It can be proven that the conditional distribution of $q$ given $(X_1 =
x_1, X_2 = x_2, \dots, X_n = x_n)$ is a beta with left parameter
$\alpha + k$ and right parameter $\beta + n - k$, with $k =
\sum_{i=1}^n x_i$. Since the prior $p(q)$ and the posterior $p(q \mid X=x)$ belong to the same
probability distribution family they are called conjugate
distributions. In particular, the beta distribution is the conjugate prior to the
Bernoulli distribution. This posterior distribution could then be used
as the new prior, without the need for expensive MCMC. The
parameters $\alpha$ and $\beta$ are updated incorporating new
experience as it comes in a computationally efficient way.

\subsection{Risk measures}
\label{sec:risk}
As we showed above, computing the posterior distribution of the continuous random
variable $Q$ is computationally efficient. A common choice for the
best estimate of the transition probability would be the maximum a
posteriori (MAP) estimate. The MAP estimate does not account for rare events and
the associated risk, an important aspect in engineering applications
where we need to be aware of the risks connected to failures. To overcome these limitations we use the conditional value at risk (CVaR) for the posterior estimate of the transition probabilities, thus enhancing the reliability and the risk-awareness of the digital twin. CVaR was first introduced for portfolio risk analysis~\cite{rockafellar2000optimization,
  rockafellar2002conditional, rockafellar2013superquantiles}, and successfully used over the years in engineering applications such as in~\cite{zhang2016decomposition, royset2017risk, yang2017algorithms, chaudhuri2022certifiable, airaudo2024risk}.

To define the CVaR we first introduce the value at risk (VaR)~\cite{pritsker1997evaluating,
  holton2003value}, also called $\alpha$-quantile. For continuous
distributions the VaR reads as follows: 
\begin{equation}
\text{VaR}_{1 - \alpha} (Q) := \inf_{t \in \mathbb{R}} \{ t :
  \text{Pr}(Q \leq t) \geq 1 - \alpha \}, \qquad \alpha \in (0, 1].
\end{equation}

CVaR, also called average value at risk or superquantile, is a coherent
risk measure, even though it is derived from the value at risk, which is
not. It measures the expected loss in the lower tail
given that a particular threshold has been crossed. It provides a
measure of the average loss magnitude beyond the $\alpha$-quantile level. 
For continuous distributions, the CVaR is defined as
\begin{equation*}
\text{CVaR}_{1 - \alpha} (Q) := \frac{1}{\alpha} \int_0^\alpha
\text{VaR}_{1-t} (Q) dt = \mathbb{E} [Q \mid Q \geq \text{VaR}_{1-\alpha}(Q)], \qquad \alpha \in (0, 1].
\end{equation*}
It follows that for every $\alpha \in (0, 1]$ this inequality
holds: $\text{VaR}_{1 - \alpha} (Q) \leq \text{CVaR}_{1 - \alpha} (Q)$.

\begin{figure}[ht]
  \centering
  \includegraphics[width=.6\textwidth, trim=10 10 10 0, clip]{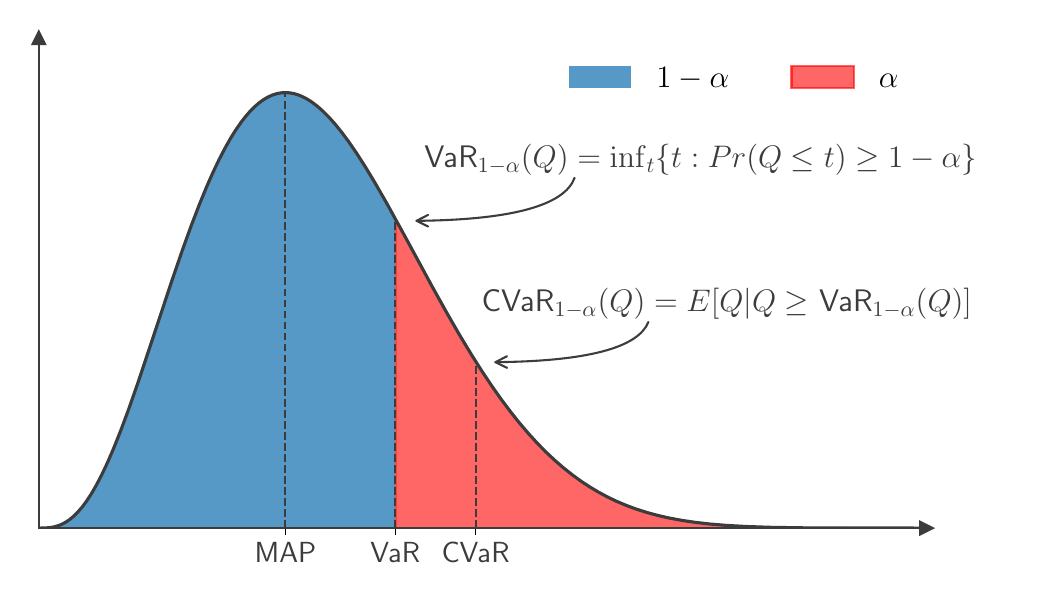}
  \caption{Maximum a posteriori (MAP) estimate, value at risk (VaR), and
    conditional value at risk (CVaR) for a generic probability 
    distribution followed by the random variable $Q$.}
  \label{fig:risk_measures}
\end{figure}

Figure~\ref{fig:risk_measures} depicts the relation between the MAP estimate, the VaR, and the CVaR for a beta-distributed random variable.

\section{Numerical results}
\label{sec:results}

This section presents two approaches to planning using robust policy optimization for a parametric MDP. The first is a risk-averse planning structure, whereby the model checking tool will compute the set of policies that are guaranteed to satisfy a given safety specification and then select the optimal. The second is used when we cannot make an almost-sure guarantee on safety but we seek to minimize the total cost of a high-impact safety incident (such as wing failure). In the next section we will describe the conditions that allow for the first form and planning and the subsequent modifications that inform the second. The probabilistic graphical model framework is implemented using the open source Python library pgmpy~\cite{pgmpy}.

\subsection{Example scenarios}

\begin{figure}[ht]
    \centering
    \subfloat[The UAV departs from a
      warehouse or a vertiport and has to deliver the package to one
      of the clients in the served area. The agent computes the best
      path and the associated maneuvers to take at every time step.
    \label{fig:grid}]{
    \includegraphics[width=.475\textwidth, trim=0 -5 0 0, clip]{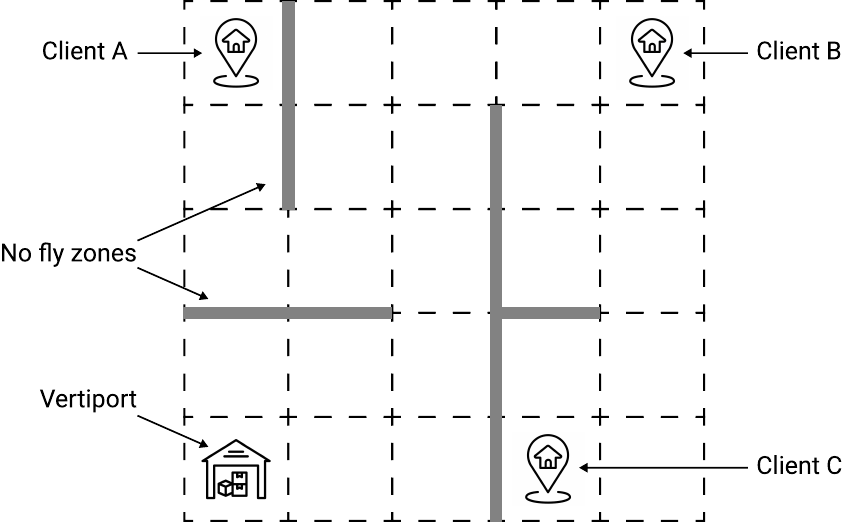}}
    \quad
   \subfloat[Aircraft attempts to
      avoid an incoming vehicle using the same $x$-$y$ trajectory. It may select a band of altitude with a gentle maneuver
      ($g_{\text{up}},~g_{\text{flat}},~g_{\text{down}}$) or an aggressive one
      ($a_{\text{up}},~a_{\text{down}}$). 
    \label{fig:vertical_air}]{
    \includegraphics[width=.475\textwidth, trim=0 25 0 0, clip]{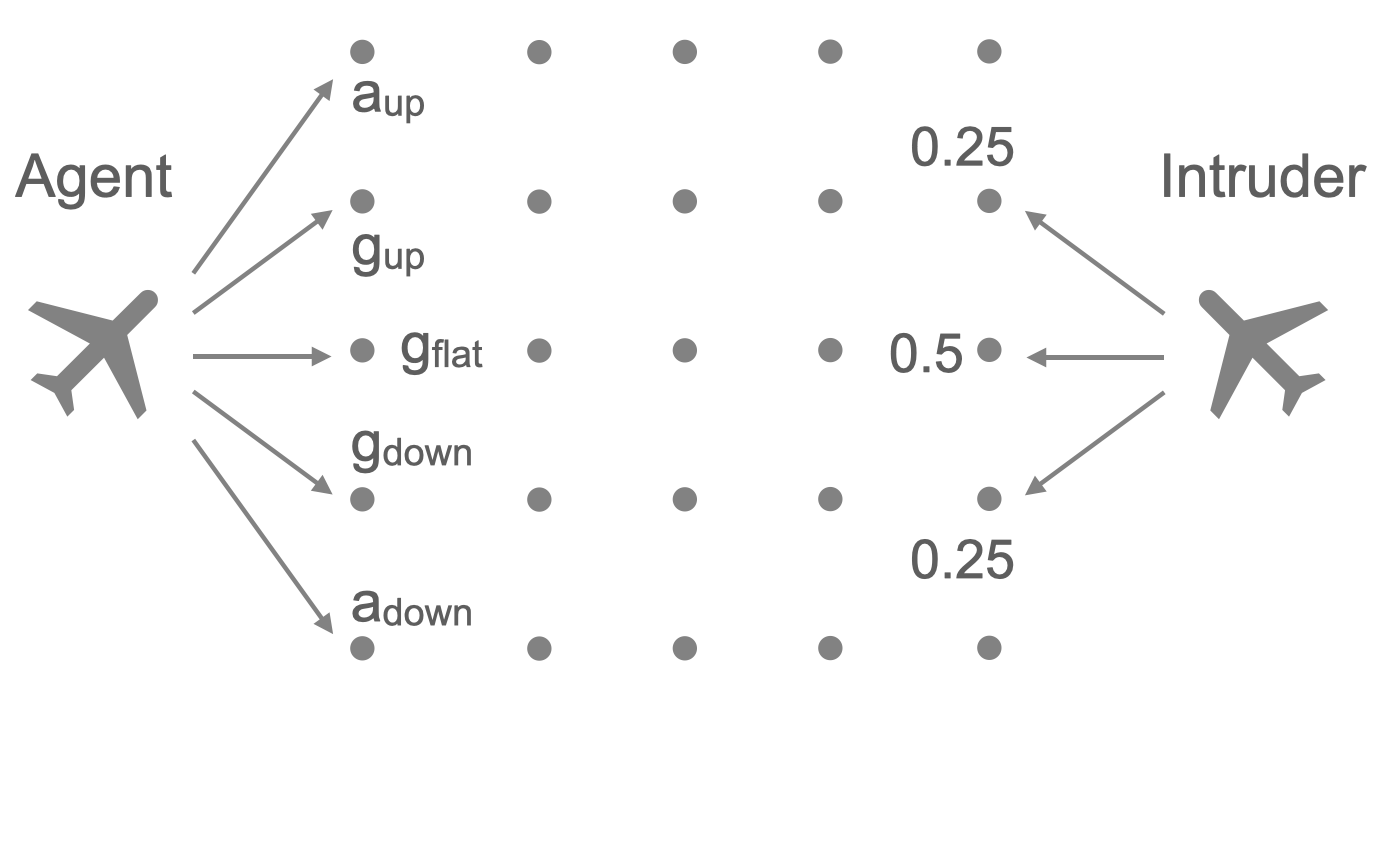}}
    \caption{(a) Last mile delivery and (b) Vertical air collision avoidance examples. }
\end{figure}

We present two case studies for UAV applications where our proposed framework provides impact when planning. These cases are described as:
\begin{enumerate}
    \item \textbf{UAV delivery} (Figure~\ref{fig:grid}) -- a UAV is delivering items in a gridworld environment with a series of target points. At each time interval there is a time-based cost for transit and a non-zero probability $q \in (0,1]$ that the wings may incur some damage in the panels, i.e., local stiffness degradation, shown in Figure~\ref{fig:pipeline_scheme}. The UAV may elect to take an \emph{aggressive} or \emph{gentle} maneuver. An aggressive maneuver results in a faster transition between states (i.e., a higher reward $r_{\textsc{agg}} > r_{\textsc{gen}}$) but it has a higher probability of incurring damage on one of the wing panels (i.e., $q_{\textsc{agg}} > q_{\textsc{gen}}$). 
    \item \textbf{Aircraft collision avoidance} (Figure~\ref{fig:vertical_air}) -- a modified MDP model~\cite{kochenderfer2012next} in a similar environment to the UAV delivery, an aircraft is attempting to deconflict with an oncoming vehicle by selecting the appropriate altitude band. At each time-interval they may select from one of five decisions for their next altitude band, three decisions $\{ g_{\text{up}},g_{\text{flat}},g_{\text{down}}\}$ form the gentle band (which carry probability $q_{\textsc{gen}}$ of causing damage to the wing panel) and two decisions $\{ a_{\text{up}},a_{\text{down}}\}$ are more aggressive maneuvers that reach further states at the cost of higher probability $q_{\textsc{agg}}$ of causing damage to the wings. The scenario is concluded after the vehicles occupy the middle band of points. If both vehicles are at the same altitude in this trajectory then they are determined to have a risk of crash and are considered safe otherwise. The opponent moves stochastically between altitude bands according to some known distribution~\cite{kochenderfer2012next}.
\end{enumerate}

These scenarios can be modeled by a parametric MDP with two parameters: the probability that the aggressive maneuver causes damage ($q_{\textsc{agg}}$) and the probability that the gentle maneuver causes damage ($q_{\textsc{gen}}$). 
There are numerous methods for synthesizing policies that ensure an agent modeled as a parametric MDP will adhere to a temporal logic specification~\cite{dehnert2015prophesy,cubuktepe2018synthesis,hahn2011synthesis}. In this work, we must both ensure a probability of reaching a
goal state with the minimization of an expected cost and therefore we employ the probabilistic model checker Storm~\cite{DehnertJK017}.
Using Storm, we instantiate the worst case values for the damage probability and formulate the problem as a linear program~\cite{BK08}. The computed policy is the policy that minimizes the expected cost to reach the goal while ensuring the UAV will remain below the damage threshold. Additionally, since Storm is an efficient probabilistic model checker, we may update the expected value of the parameters as the aircraft makes decisions and computes updates to the policy online.

In these examples, the worst case value of the parameter is
predictable and therefore we chose an efficient tool to compute the
policy. When it is uncertain how a given parameter value will affect
the property, one may employ methods such as sequential convex
programming~\cite{CubuktepeJJKT22}. For other approaches for verification of uncertain parametric MDPs see~\cite{badings2022scenario}.

In each case above the objective of the aircraft can be described using a temporal logic formula such as $\phi = \neg C \cup G $, i.e., do not crash until you reach the goal where $C$ describes a crash state such as wing damage, occupying the same altitude band as an opponent or violating a no-fly region, and $G$ is some goal condition such as delivering an object or passing an opponent in a different altitude band.

In both scenarios we seek to ensure that: 
\begin{enumerate}
    \item we successfully navigate to a goal $G$;
    \item  without incurring catastrophic structural damage to our wings $\neg C$;
    \item with the minimum expected cost of arrival.
\end{enumerate}
At each decision point, there is a trade-off between taking a cheap, aggressive but higher-risk action against an expensive, gentle but lower-risk action. For an almost-sure guarantee on arriving to the goal without incurring catastrophic damage there needs to be either a zero probability action (i.e., $q_{\textsc{gen}}=0$) for wing damage or a minimum number of decisions smaller than the number of damage bins on the wing (e.g., 5 steps to the goal and 10 bins for structural damage transitions to occur). While there exist many situations where $q_{\textsc{gen}}=0$, we are most interested in presenting the trade-off between ensuring safety and optimizing for time taken as well as showing the value of the predictive input of an updating digital twin on the plan. Accordingly, most results, unless otherwise specified will be minimizing the expected cost, where the cost of a failure condition is graded as significantly higher than any time-based action.

For all the test cases we assume that the calibration of the structural DT has already been performed as in~\cite{kapteyn2021probabilistic}. The initial prior on the transition probabilities is provided as a beta distribution $\mathcal{B}e(2, \beta)$ with left parameter $\alpha = 2$ and given mode. 
The formula to compute the mode of a beta-distributed random variable given its parameters is the following:
\begin{equation}
\label{eq:betamode}
    \frac{\alpha - 1}{\alpha + \beta - 2}, \qquad \text{for } \, \alpha, \beta > 1.
\end{equation}
By using Equation~\eqref{eq:betamode} we compute the corresponding right parameter $\beta$. To ensure the computational model interpretability we round $\beta$ to the closest integer.

\subsection{Digital state estimation}
\label{sec:inverse}
This section describes how the digital state is estimated from the noisy observational data coming from the sensors placed on the wing of the UAV. We emphasize there are many different methods that can be leveraged to solve the inverse problem from data to parameters comprising the digital state. The reader can refer to Bayesian state estimation techniques or more generally to data assimilation methods~\cite{wikle2007bayesian, reich2015probabilistic}. The crucial aspect is to compute the uncertainty associated with the estimations to fully leverage the graphical model.
In this application the digital state is a four-dimensional vector composed by $\mathbf{x} \in \mathbb{N}^2$, representing the deterministic UAV position coordinates, and by $\mathbf{z} \in \mathbb{D} \subset \mathbb{R}^2$, which is a vector comprising the structural health parameters associated to the UAV wing. Since the UAV position is not affected by uncertainty in our setting, we focus on $\mathbf{z}$. The vector $\mathbf{z}$ is defined as
\begin{equation}
\mathbf{z} = [z_1, z_2] \in \mathbb{D} \subset \mathbb{R}^2, \qquad \mathbb{D} := \left \{ \left (\frac{i}{10}, \frac{j}{10} \right) \bigg | \, i, j \in [0, 1, \dots, 8] \right \}.
\end{equation}
It represents the percentage of stiffness degradation (from 0\% to 80\%) of the two red regions on the wing showed in Figure~\ref{fig:pipeline_scheme}.
The vector $\mathbf{z}$ is estimated at every time step through a set of strain measurements $\boldsymbol{\epsilon} = \{ \epsilon^j \}_{j=1}^{24}$, coming from sensors placed on the UAV. We assume these measurements are corrupted with additive white Gaussian noise:
\begin{equation}
\hat{\epsilon}^j_t = \epsilon^j_t + v_t, \qquad v_t \sim
\mathcal{N}(0, \sigma),
\end{equation}
where the subscript denotes the time instant, and $\sigma=10$ so that the sensor model does not match the simulated measurements.
To estimate $\mathbf{z}$ we solve the following inverse problem
\begin{equation}
\min_{\boldsymbol{\theta}} \frac{1}{2} \| F(\boldsymbol{\theta}) - \hat{\boldsymbol{\epsilon}} \|^2_2 + \| \boldsymbol{\theta} \|_2,
\end{equation}
where $\boldsymbol{\theta} \in [0, 0.8]^2$ and $\hat{\boldsymbol{\epsilon}} \in \mathbb{R}^{24}$ is the vector comprising the noisy strain measurements. We use an $L_2$ regularization to favor low norm solutions. The $\theta^*$ that realizes the minimum is then approximated with the nearest point in the discrete space $\mathbb{D}$.
The map $F: [0, 0.8]^2 \to \mathbb{R}^{24}$ describes the evaluation
of the reduced order model based on the static condensation reduced
basis element (SCRBE) method developed in~\cite{huynh2013static, eftang2013port, smetana2016optimal} and presented in~\cite{kapteyn2020toward, kapteyn2020data} for the UAV we are considering. The SCRBE technique brings together the scalability of component-based formulations and the accuracy of the reduced basis method. 

The resulting accuracy of the state estimation is equal to $75.4 \%$. The accuracy is measured by sampling $100$ realizations of the noisy sensor data and then computing the average of the confusion matrix' diagonal. In particular, if we look at the estimation of the single components of the vector $\mathbf{z}$, we observe that $z_1$ is estimated with a $100\%$ accuracy, and all the errors are due to $z_2$, which is in line with the findings in~\cite{kapteyn2021probabilistic}. These measures are factored in the PGM through $\phi_t^{\text{data}}$.

\subsection{Risk-averse case using CVaR}
\begin{figure}[ht]
    \centering
    \includegraphics[width=.75\textwidth, trim=18 0 0 0, clip]{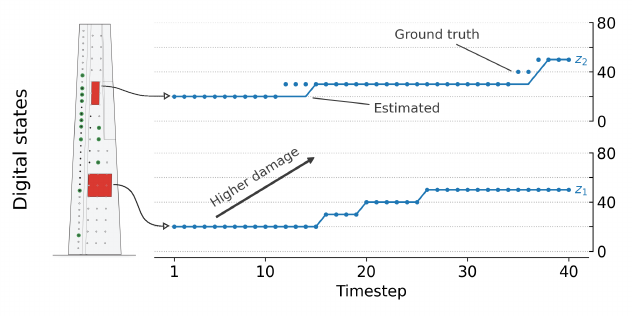}
    \caption{Evolution of the digital states $z_1$ and $z_2$ for the
      test mission. The initial state $\mathbf{z}$ at time $t=1$ is set to $[0.2, 0.2]$.}
    \label{fig:state_evolution}
\end{figure}

 The simulated mission is composed of $40$ time steps and it starts from an initial state $\mathbf{z} = [0.2, 0.2]$. In Figure~\ref{fig:state_evolution} we show how the data assimilation performs for $z_1$ (bottom part) and $z_2$ (top part). The underlying transition probabilities are $3\%$ and $10\%$ for the gentle and aggressive actions, respectively. As we can see the tracking of $z_1$ is perfect, while for $z_2$ we need some time steps to identify a state change. This is observed also in~\cite{kapteyn2021probabilistic} and it is due to the chosen state discretization levels and the position on the wing, closer to the tip, which is more sensitive to noise.

\begin{figure}[ht]
    \centering
    \includegraphics[width=.75\textwidth]{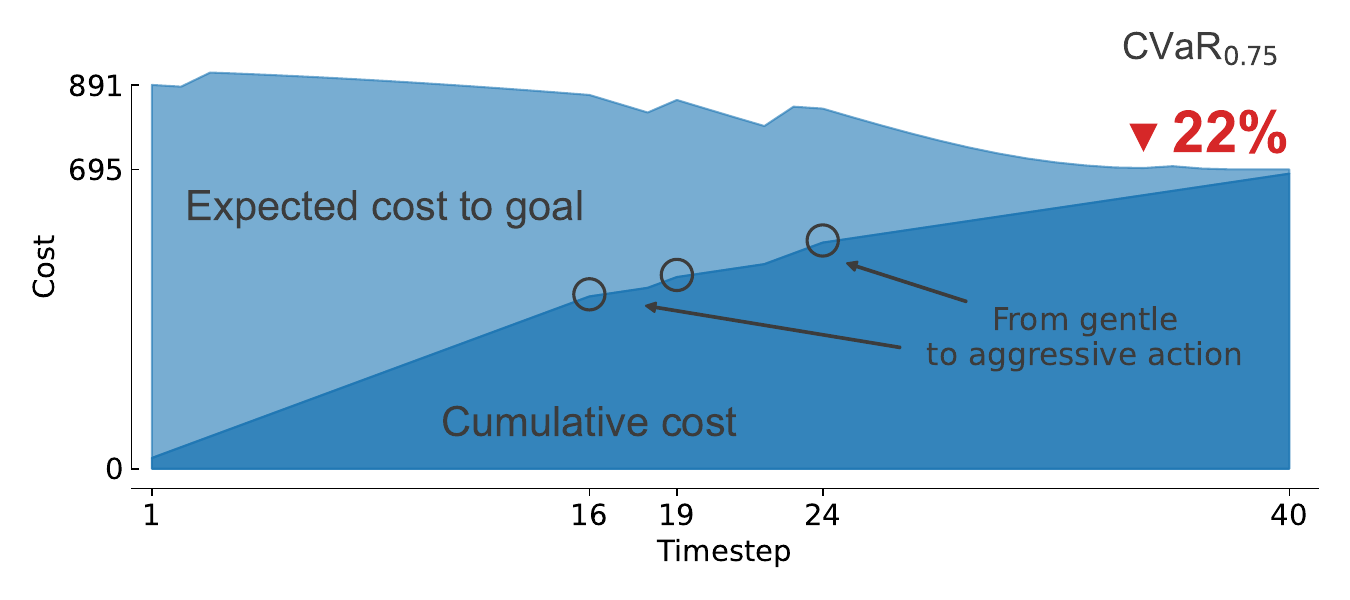}
    \caption{Expected and effective cost to complete a mission using the conditional value at risk to estimate the transition probability matrices.}
    \label{fig:mission_cost_q}
  \end{figure}

\begin{figure}[htb]
    \centering
    \includegraphics[width=.49\textwidth]{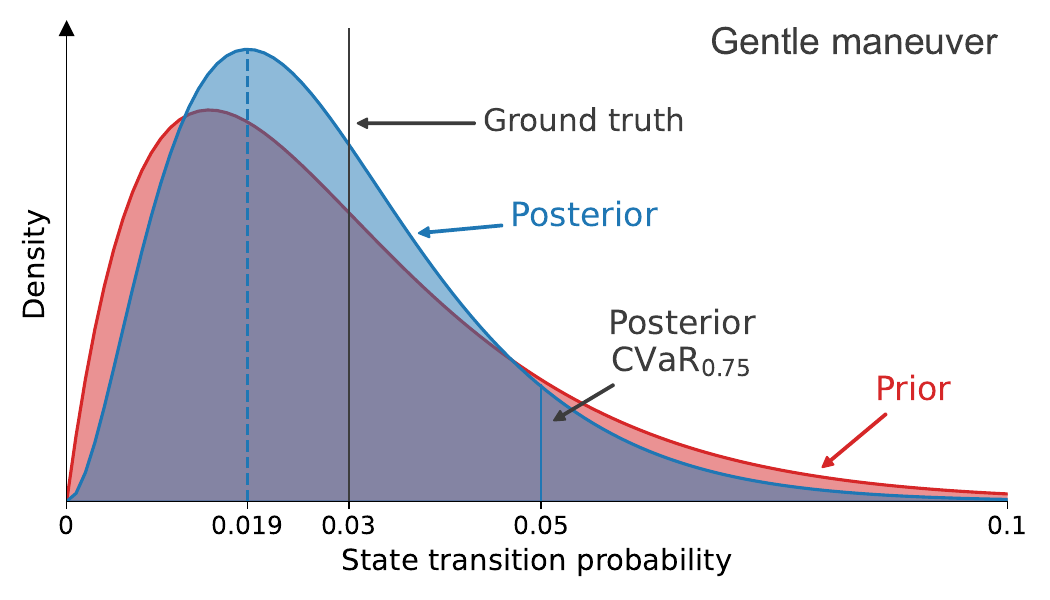}\hfill
    \includegraphics[width=.49\textwidth]{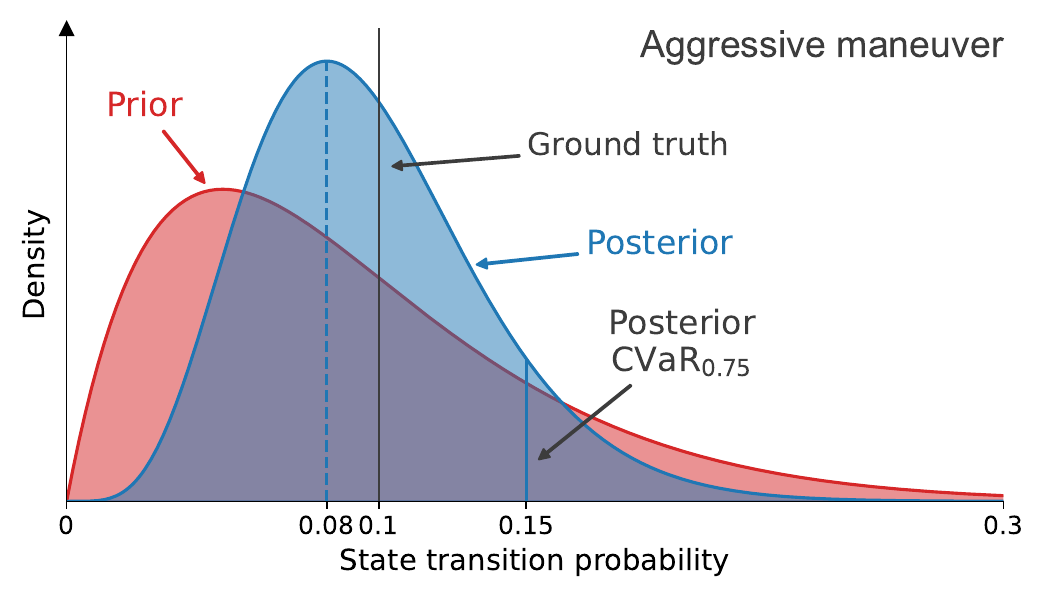}
    \caption{Prior and posterior distributions corresponding to the gentle maneuver (left panel) and to the aggressive maneuver (right panel) for the risk-averse test case.}
    \label{fig:betas_q}
\end{figure}

In Figure~\ref{fig:mission_cost_q} we report the expected cost to complete the mission at every time step, and the actual cumulative cost spent to follow the optimal policy. To estimate the transition probabilities we use the conditional value at risk with $\alpha=0.25$, i.e., the CVaR at the upper quantile. The adaptive policy results in a reduced cost equal to $22\%$. This reduction is due to the policy refinement and the posterior updates. If we consider the last-mile delivery scenario, the cost reduction translates to faster deliveries without incurring a fatal crash, while accounting for the risk associated with extreme events. By looking at the actual cost slope we can identify the switch between gentle and aggressive actions, which correspond to a cost equal to 25 and 10, respectively. The priors associated to the gentle and aggressive maneuvers are $\mathcal{B}e(2, 66)$ and $\mathcal{B}e(2, 20)$, respectively. They are reported in Figure~\ref{fig:betas_q} together with the corresponding posteriors after 40 time steps. We can see that in both cases we are converging to the unknown underlying transitions.

\subsection{Best estimate case using MAP}
For this test case we use the maximum a posteriori estimate of the beta distributions to approximate the transition probabilities. 
The initial state at $t=1$ and the costs associated with the two actions are the same as in the previous test case. The underlying transition probability for the gentle maneuver is equal to $2\%$, while for the aggressive maneuver is $10\%$. 

\begin{figure}[ht]
    \centering
    \includegraphics[width=.75\textwidth]{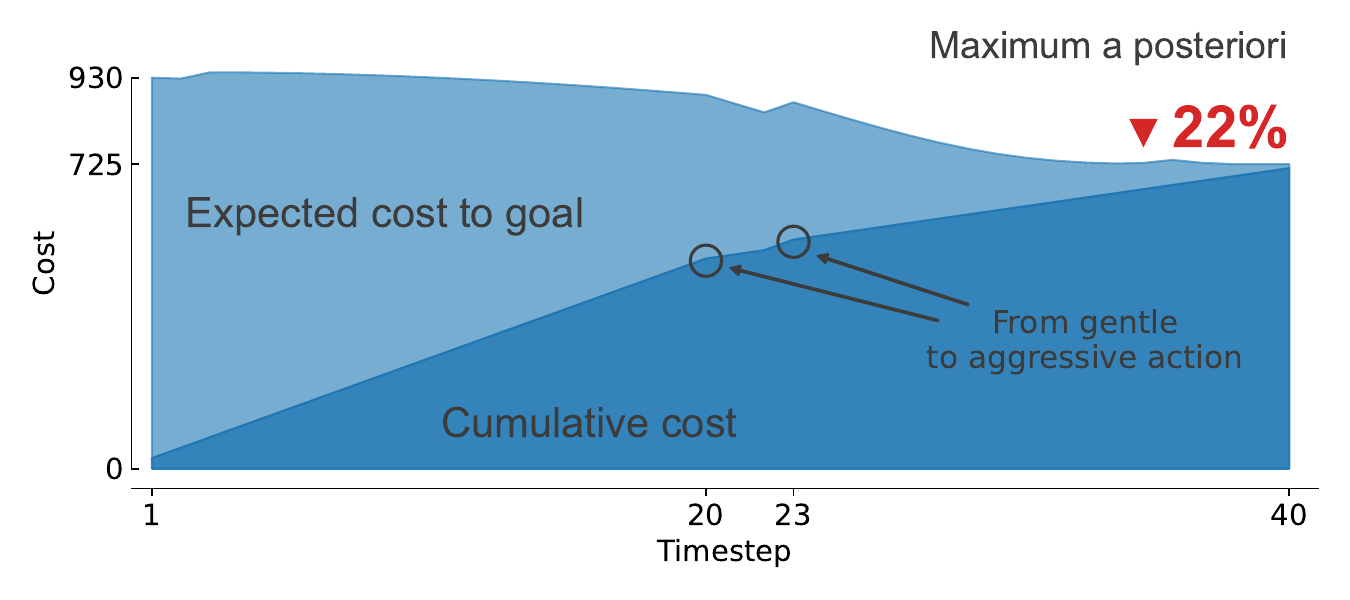}
    \caption{Expected and effective cost to complete a mission using
      the maximum a posteriori estimate to estimate the transition probability matrices.}
    \label{fig:mission_cost_map}
\end{figure}

By comparing the expected cost to goal and the actual cumulative cost we notice a similar reduction as before equal to $22.07\%$, as we can see in Figure~\ref{fig:mission_cost_map}. The similarities are due to the fact that we have a limited set of possible actions and the length of the mission is also the same in the two scenarios. 

\begin{figure}[htb]
    \centering
    \includegraphics[width=.49\textwidth]{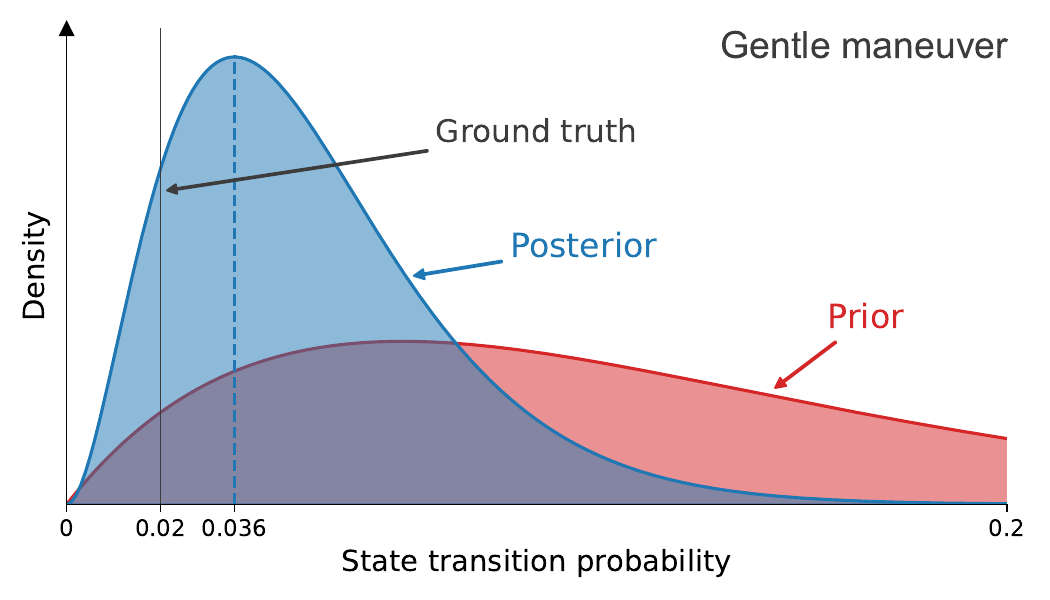}\hfill
    \includegraphics[width=.49\textwidth]{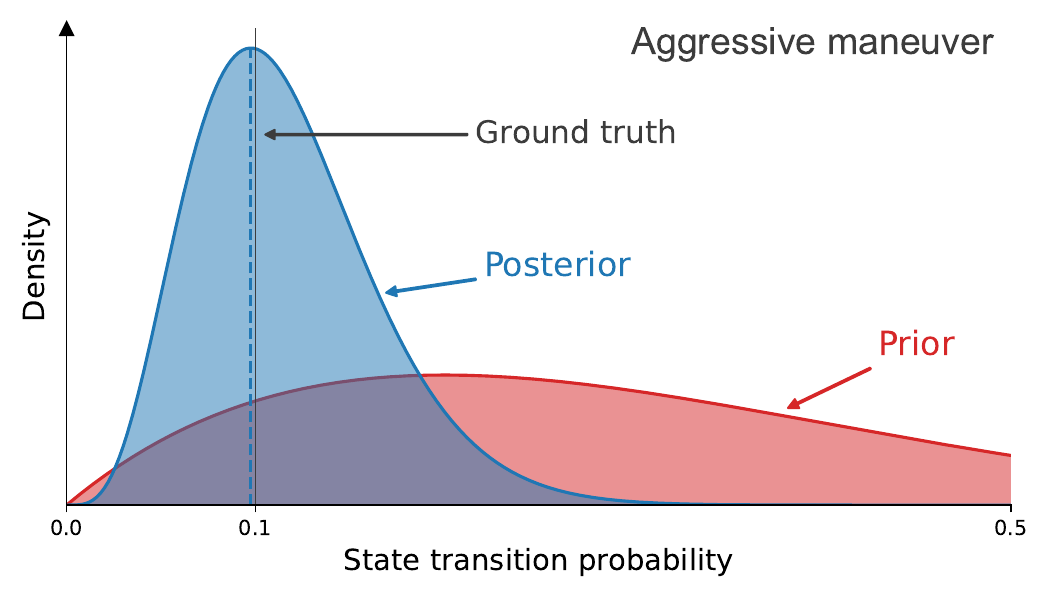}
    \caption{Prior and posterior distributions corresponding to the
      gentle maneuver (left panel) and to the aggressive maneuver
      (right panel) for the maximum a posteriori estimate test case.}
    \label{fig:betas_map}
\end{figure}

The priors for gentle and aggressive actions are $\mathcal{B}e(2, 14)$ and $\mathcal{B}e(2, 5)$, respectively. They correspond to a MAP equal to $7\%$ and $20\%$. In Figure~\ref{fig:betas_map} we also plotted the posteriors after the last time step $t=40$. For the aggressive maneuver the MAP estimate identifies almost perfectly the ground truth mode, while for the gentle action we need more time steps.

\subsection{Digital state predictions}
\begin{figure}[ht]
    \centering
    \includegraphics[width=.325\textwidth, trim=75 0 30 5, clip]{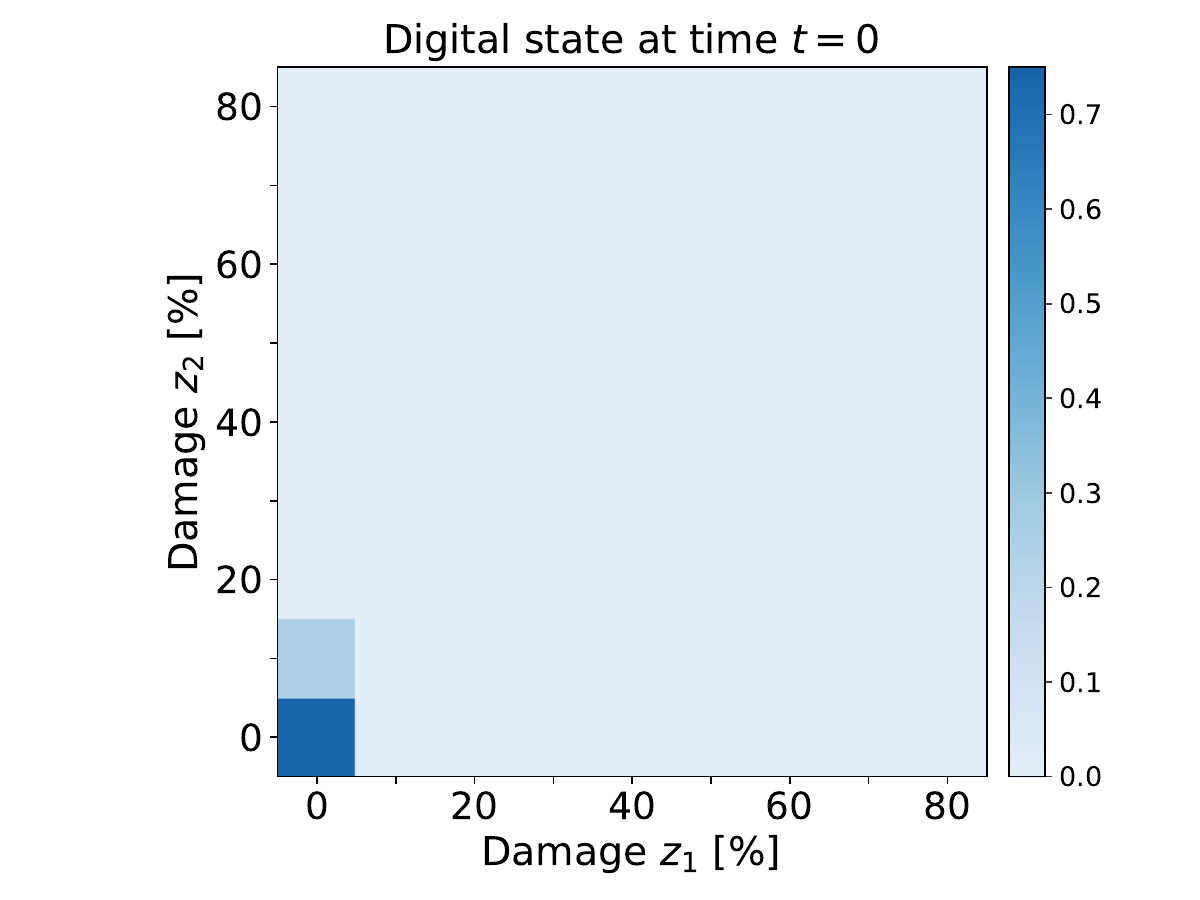}
    \includegraphics[width=.325\textwidth, trim=75 0 30 5, clip]{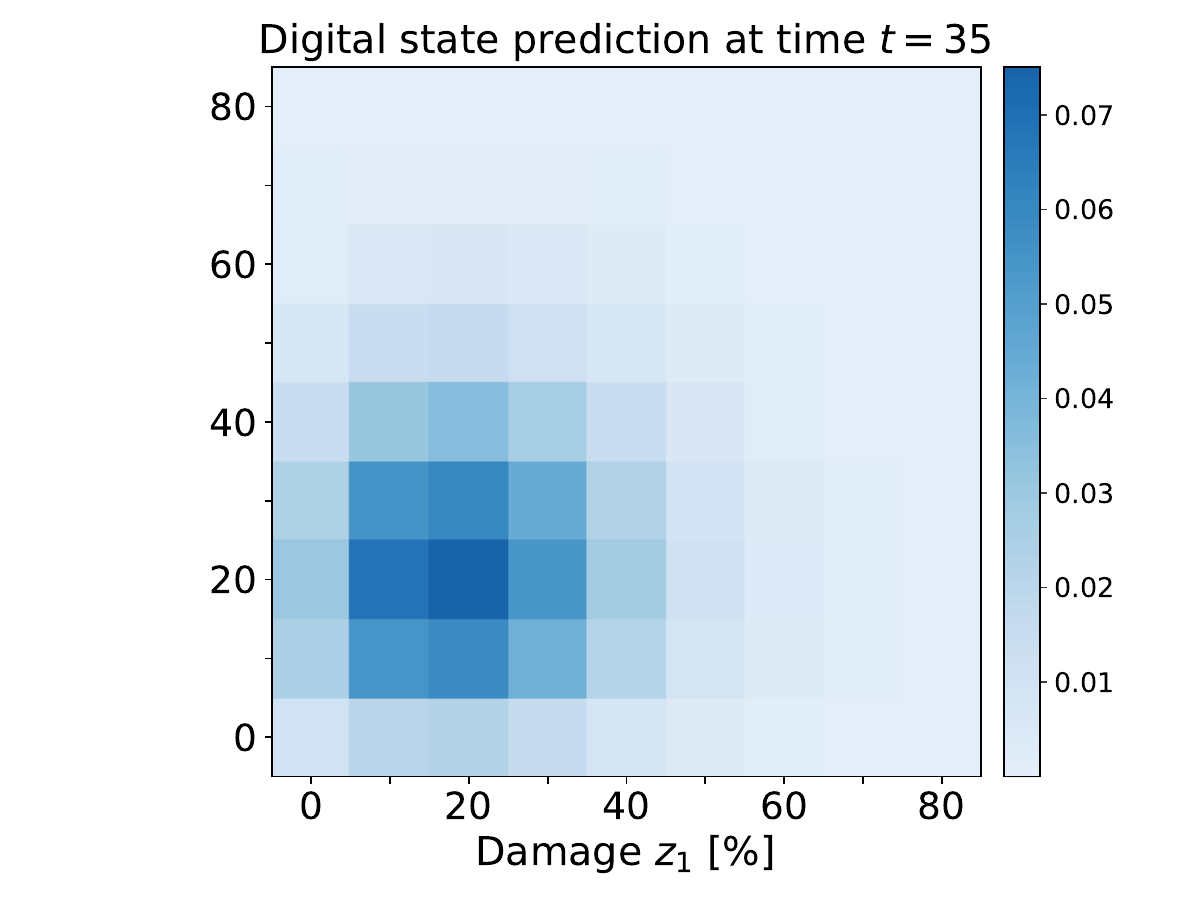}
    \includegraphics[width=.325\textwidth, trim=75 0 30 5, clip]{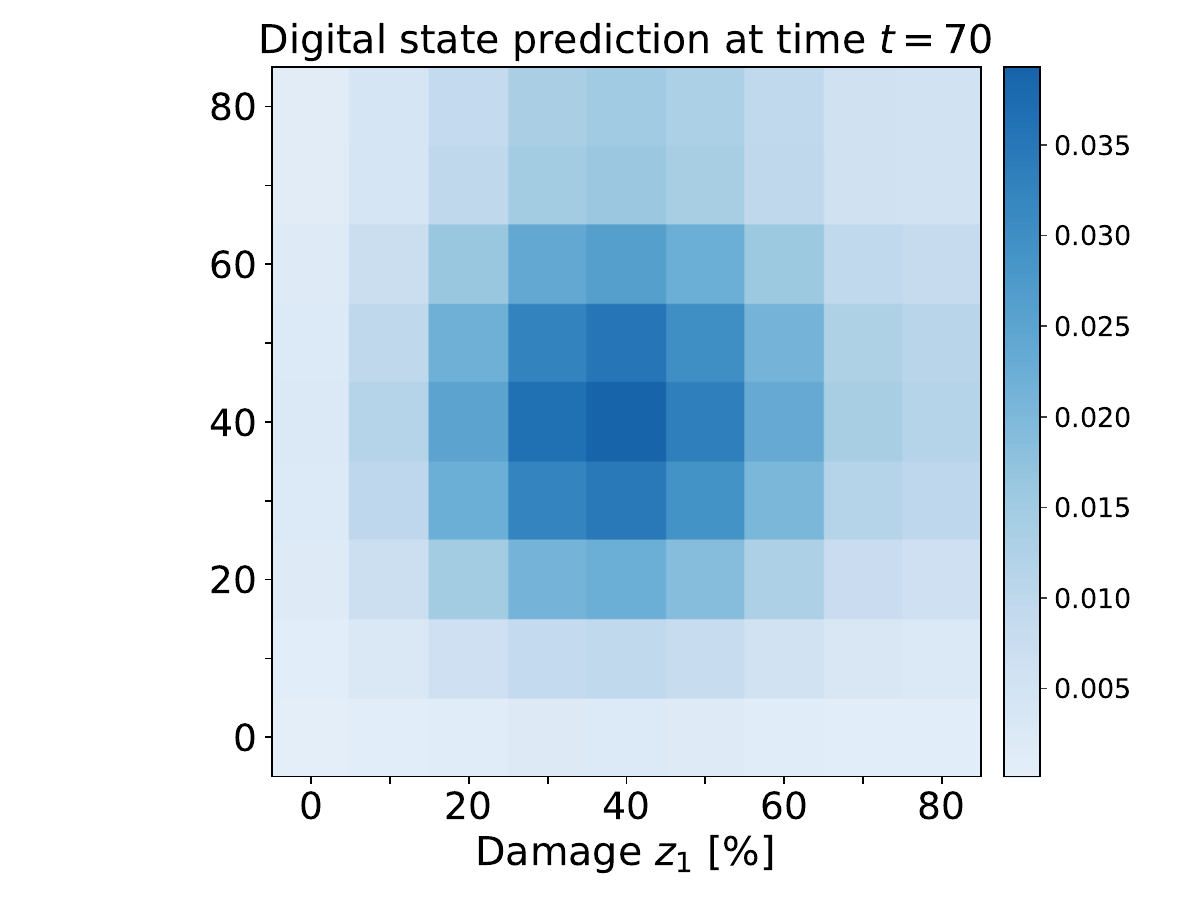}
    \caption{Digital state predictions. Initial estimate in the left panel ($t=0$), predictions after $35$ time steps in the central panel, and after $70$ time steps in the right panel. The colorbar represents the probability of being in a particular state $[z_1, z_2]$.}
    \label{fig:predictions}
\end{figure}

The PGM allows making predictions by propagating all the uncertainties in the graph. By leveraging the subgraph presented in Figure~\ref{fig:prediction_graph}, we compute the evolution of the digital state with the policy $\pi_{t_c}$.  Figure~\ref{fig:predictions} depicts the digital state starting from $t_c=0$, where we have a discrete distribution corresponding to a $75\%$ probability of being at the state $[0, 0]$ and a $25\%$ probability of being in $[0, 0.1]$. The prediction horizon is extended over $70$ time steps in the future so that $t_p=t_c+70$. In the central panel we show the estimates at half of the time horizon, while in the right panel the final estimates. In this example we employ a posterior estimate of the transition probabilities equal to $2\%$ and $10\%$ for the gentle and aggressive actions, respectively. The DT formulation informs about the expected future degradation of structural health, allowing to planning of preventive interventions and maintenance.

\section{Conclusions and future perspectives}
\label{sec:the_end}
In this paper we integrated uncertainty quantification, risk-awareness, adaptive planning for autonomous systems, and predictive capabilities using a probabilistic graphical model for digital twins. We showed how it is possible to dynamically update the underlying assumptions in the creation of a PGM using parametric Markov decision processes. We demonstrated that the decision-making agent will deploy a robust policy while dynamically refining its belief from new sensor data. The dynamic update of the policy is a crucial characteristic that every digital twin should have since it increases personalization and provides better decision-making while accounting for risk and rare events. Not only is the DT able to monitor the physical asset in a more realistic way, but also the graphical model improves the DT’s prediction accuracy. The increased personalization of the DT represents a stepping stone for enhanced predictive maintenance and monitoring. 

The need of recomputing the optimal policy at every time step could be
computationally demanding depending on the dimension of the discretized variables and
on the actual computational power on board. Therefore, the ability to
accurately describe the real-world scenario by the MDP is limited by
the actual controller installed on the UAV.
Nevertheless, the methods presented in this work aim at setting a general
formulation for digital twins that are able to adapt to sensed scenarios
by recomputing optimal policies, and accounting for risk. Future studies will be conducted considering more general parametrized transition probabilities.

\section*{Acknowledgements}
The authors acknowledge support from the NASA University Leadership
Initiative under Cooperative Agreement 80NSSC21M0071. The authors also
thank Michael Kapteyn for the UAV's wing illustrations, and Akselos
for the structural analysis and reduced order modeling software used
(Akselos Integra v4.5.9).

\bibliographystyle{abbrvurl}

\end{document}

%% file: figures/digital_twin_graph.tex
\tikzstyle{P_node} = [circle,draw=black,thick,align=center,minimum size=0.8cm]
\tikzstyle{O_node} = [circle,draw=black,thick, line width=.55mm, align=center,minimum size=0.8cm]
\tikzstyle{D_node} = [circle,draw=black,thick,align=center,minimum size=0.8cm]
\tikzstyle{D_A_node} = [circle,draw=black,thick,align=center,minimum size=0.8cm]
\tikzstyle{U_node_prob} = [circle, draw=black, thick, minimum size=0.8cm]
\tikzstyle{U_node_act} = [rectangle, draw=black, thick, line width=.55mm, minimum width=0.8cm, minimum height=0.8cm]
\tikzstyle{R_node} = [rectangle,draw=black,thick,align=center,minimum width=0.8cm, minimum height=0.8cm,rotate=45]

\begin{tikzpicture}[scale=.85, every node/.style={scale=.9}]

\node [rotate=90] () at (-1.5, 4.25) {Physical Space};
\node [rotate=90] () at (-1.5,.8) {Digital Space};

\node [D_node] (D_0) at (0,0.3) {};
\node [] at (0,0.3) {$D_{0}$};
\node [R_node] (R_0) at (2,1.3) {};
\node [] at (2,1.3) {$R_{0}$};
\node [P_node] (S_0) at (0,5) {};
\node [] at (0,5) {$S_{0}$};
\node [O_node] (O_0) at (0,3.5) {};
\node [] at (0,3.5) {$O_{0}$};
\node [U_node_act] (U_0) at (2,3.5) {};
\node [] at (2,3.5) {$U_{0}$};

\node [D_node] (D_1) at (4,0.3) {};
\node [] at (4,0.3) {$D_{1}$};
\node [R_node] (R_1) at (6,1.3) {};
\node [] at (6,1.3) {$R_{1}$};
\node [P_node] (S_1) at (4,5) {};
\node [] at (4,5) {$S_{1}$};
\node [O_node] (O_1) at (4,3.5) {};
\node [] at (4,3.5) {$O_{1}$};
\node [U_node_act] (U_1) at (6,3.5) {};
\node [] at (6,3.5) {$U_{1}$};

\node [D_node] (D_2) at (8,0.3) {};
\node [] at (8,0.3) {$D_{2}$};
\node [R_node] (R_2) at (10,1.3) {};
\node [] at (10,1.3) {$R_{2}$};
\node [P_node] (S_2) at (8,5) {};
\node [] at (8,5) {$S_{2}$};
\node [O_node] (O_2) at (8,3.5) {};
\node [] at (8,3.5) {$O_{2}$};
\node [U_node_act] (U_2) at (10,3.5) {};
\node [] at (10,3.5) {$U_{2}$};

\node [] (midstart) at (-2,2.5) {};
\node [] (midfinal) at (12,2.5) {};

\node [] (t0) at (-0.8,-1.2) {};
\node [] () at (0,-1.2) {$|$};
\node [] () at (0,-0.7) {$t = 0$};
\node [] () at (4,-1.2) {$|$};
\node [] () at (4,-0.7) {$t = 1$};
\node [] () at (8,-1.2) {$|$};
\node [] () at (8,-0.7) {$t = 2$};
\node [] (t2) at (12,-1.2) {};

\draw[-latex,thick,black] (S_0) to (O_0);
\draw[-latex,thick,black] (O_0) to (D_0);
\draw[-latex,thick,black] (D_0) to (U_0);
\draw[-latex,thick,black] (D_0) to (R_0);
\draw[-latex,thick,black] (U_0) to (R_0);
\draw[-latex,thick,black] (U_0) to (S_1);
\draw[-latex,thick,black] (U_0) to (D_1);
\draw[-latex,thick,black] (D_0) to (D_1);
\draw[-latex,thick,black] (S_0) to (S_1);

\draw[-latex,thick,black] (S_1) to (O_1);
\draw[-latex,thick,black] (O_1) to (D_1);
\draw[-latex,thick,black] (D_1) to (U_1);
\draw[-latex,thick,black] (D_1) to (R_1);
\draw[-latex,thick,black] (U_1) to (R_1);
\draw[-latex,thick,black] (U_1) to (S_2);
\draw[-latex,thick,black] (U_1) to (D_2);
\draw[-latex,thick,black] (D_1) to (D_2);
\draw[-latex,thick,black] (S_1) to (S_2);

\draw[-latex,thick,black] (S_2) to (O_2);
\draw[-latex,thick,black] (O_2) to (D_2);
\draw[-latex,thick,black] (D_2) to (U_2);
\draw[-latex,thick,black] (D_2) to (R_2);
\draw[-latex,thick,black] (U_2) to (R_2);

\draw[thick,dotted,black] (S_2) to (11.9, 5);
\draw[thick,dotted,black] (U_2) to (11.85, 0.87);
\draw[thick,dotted,black] (D_2) to (11.9, 0.3);

\draw[black!70] (midstart) to (midfinal);
\draw[-latex,thick,black] (t0) to (t2);

\end{tikzpicture}

%% file: figures/prediction_graph.tex
\tikzstyle{D_node} = [circle,draw=black,thick,align=center,minimum size=0.95cm]
\tikzstyle{U_node_prob} = [circle, draw=black, thick, minimum size=0.95cm]

\begin{tikzpicture}[scale=.9, every node/.style={scale=0.9}]

\node [rotate=90] () at (-1.8,.8) {Digital Space};

\node [D_node] (D_0) at (0,0) {};
\node [] at (0,0) {$D_{t_c}$};
\node [U_node_prob] (U_0) at (1.5,1.5) {};
\node [] at (1.5,1.5) {$U_{t_c}$};

\node [D_node] (D_1) at (3,0) {};
\node [] at (3,0) {$D_{t_c + 1}$};
\node [U_node_prob] (U_1) at (4.5,1.5) {};
\node [] at (4.5,1.5) {$U_{t_c + 1}$};

\node [D_node] (D_2) at (6,0) {};
\node [] at (6,0) {$D_{t_c + 2}$};
\node [U_node_prob] (U_2) at (7.5,1.5) {};
\node [] at (7.5,1.5) {$U_{t_c + 2}$};

\node [D_node] (D_3) at (9,0) {};
\node [] at (9,0) {$D_{t_c + 3}$};

\node [] (t0) at (-2,-1.5) {};
\node [] () at (0,-1.5) {$|$};
\node [] () at (0,-1.) {$t=t_c$};
\node [] () at (3,-1.5) {$|$};
\node [] () at (3,-1.) {$t=t_c+1$};
\node [] () at (6,-1.5) {$|$};
\node [] () at (6,-1.) {$t=t_c+2$};
\node [] () at (9,-1.5) {$|$};
\node [] () at (9,-1.) {$t=t_c+3$};
\node [] (t2) at (11,-1.5) {};

\draw[-latex,thick,dotted,black] (-1.2, 0) to (-.5,0);
\draw[-latex,thick,dotted,black] (-1.2, 1.8) to (-.35,.35);
\draw[-latex,thick,dotted,black] (0, 1.9) to (0,.5);

\draw[-latex,thick,black] (D_0) to (U_0);
\draw[-latex,thick,black] (U_0) to (D_1);
\draw[-latex,thick,black] (D_0) to (D_1);

\draw[-latex,thick,black] (D_1) to (U_1);
\draw[-latex,thick,black] (U_1) to (D_2);
\draw[-latex,thick,black] (D_1) to (D_2);

\draw[-latex,thick,black] (D_2) to (U_2);
\draw[-latex,thick,black] (U_2) to (D_3);
\draw[-latex,thick,black] (D_2) to (D_3);

\draw[-latex,thick,black] (t0) to (t2);

\end{tikzpicture}

%% file: main_with_bib.bbl
\begin{thebibliography}{10}

\bibitem{aiaa2020digital}
{AIAA Digital Engineering Integration Committee}.
\newblock Digital twin: Definition \& value.
\newblock {\em AIAA and AIA Position Paper}, 2020.

\bibitem{airaudo2024risk}
F.~Airaudo, H.~Antil, R.~Lohner, and U.~Rakhimov.
\newblock {\em On the Use of Risk Measures in Digital Twins to Identify
  Weaknesses in Structures}.
\newblock AIAA SciTech Forum. American Institute of Aeronautics and
  Astronautics, 2024.
\newblock \href {https://doi.org/10.2514/6.2024-2622}
  {\path{doi:10.2514/6.2024-2622}}.

\bibitem{pgmpy}
A.~Ankan and A.~Panda.
\newblock {pgmpy: Probabilistic graphical models using Python}.
\newblock In {\em Proceedings of the 14th Python in Science Conference (SCIPY
  2015)}. Citeseer, 2015.

\bibitem{arcieri2023bridging}
G.~Arcieri, C.~Hoelzl, O.~Schwery, D.~Straub, K.~G. Papakonstantinou, and
  E.~Chatzi.
\newblock Bridging pomdps and bayesian decision making for robust maintenance
  planning under model uncertainty: An application to railway systems.
\newblock {\em Reliability Engineering \& System Safety}, page 109496, 2023.
\newblock \href {https://doi.org/10.1016/j.ress.2023.109496}
  {\path{doi:10.1016/j.ress.2023.109496}}.

\bibitem{badings2022scenario}
T.~Badings, M.~Cubuktepe, N.~Jansen, S.~Junges, J.-P. Katoen, and U.~Topcu.
\newblock {Scenario-based verification of uncertain parametric MDPs}.
\newblock {\em International Journal on Software Tools for Technology
  Transfer}, 24(5):803--819, 2022.
\newblock \href {https://doi.org/10.1007/s10009-022-00673-z}
  {\path{doi:10.1007/s10009-022-00673-z}}.

\bibitem{BK08}
C.~Baier and J.~Katoen.
\newblock {\em {Principles of Model Checking}}.
\newblock {MIT Press}, Cambridge, MA, 2008.

\bibitem{chaudhuri2022certifiable}
A.~Chaudhuri, B.~Kramer, M.~Norton, J.~O. Royset, and K.~E. Willcox.
\newblock Certifiable risk-based engineering design optimization.
\newblock {\em AIAA Journal}, 60(2):551--565, 2022.
\newblock \href {https://doi.org/10.2514/1.J060539}
  {\path{doi:10.2514/1.J060539}}.

\bibitem{CubuktepeJJKT22}
M.~Cubuktepe, N.~Jansen, S.~Junges, J.~Katoen, and U.~Topcu.
\newblock {Convex Optimization for Parameter Synthesis in MDPs}.
\newblock {\em {IEEE Transactions Automatic Control}}, 67(12):6333--6348, 2022.
\newblock \href {https://doi.org/10.1109/TAC.2021.3133265}
  {\path{doi:10.1109/TAC.2021.3133265}}.

\bibitem{cubuktepe2018synthesis}
M.~Cubuktepe, N.~Jansen, S.~Junges, J.-P. Katoen, and U.~Topcu.
\newblock {Synthesis in pMDPs: A tale of 1001 parameters}.
\newblock In {\em {International Symposium on Automated Technology for
  Verification and Analysis}}, pages 160--176. Springer, 2018.
\newblock \href {https://doi.org/10.1007/978-3-030-01090-4\_10}
  {\path{doi:10.1007/978-3-030-01090-4\_10}}.

\bibitem{dean1989model}
T.~Dean and K.~Kanazawa.
\newblock A model for reasoning about persistence and causation.
\newblock {\em Computational Intelligence}, 5(2):142--150, 1989.
\newblock \href {https://doi.org/10.1111/j.1467-8640.1989.tb00324.x}
  {\path{doi:10.1111/j.1467-8640.1989.tb00324.x}}.

\bibitem{dehnert2015prophesy}
C.~Dehnert, S.~Junges, N.~Jansen, F.~Corzilius, M.~Volk, H.~Bruintjes, J.-P.
  Katoen, and E.~{\'A}brah{\'a}m.
\newblock {Prophesy: A probabilistic parameter synthesis tool}.
\newblock In {\em {Computer Aided Verification: 27th International Conference,
  CAV 2015, San Francisco, CA, USA, July 18-24, 2015, Proceedings, Part I 27}},
  pages 214--231. Springer, 2015.
\newblock \href {https://doi.org/10.1007/978-3-319-21690-4\_13}
  {\path{doi:10.1007/978-3-319-21690-4\_13}}.

\bibitem{DehnertJK017}
C.~Dehnert, S.~Junges, J.~Katoen, and M.~Volk.
\newblock A storm is coming: A modern probabilistic model checker.
\newblock In {\em {CAV} {(2)}}, volume 10427 of {\em Lecture Notes in Computer
  Science}, pages 592--600, Cham, 2017. Springer.
\newblock \href {https://doi.org/10.1007/978-3-319-63390-9\_31}
  {\path{doi:10.1007/978-3-319-63390-9\_31}}.

\bibitem{delahaye2011decision}
B.~Delahaye, K.~G. Larsen, A.~Legay, M.~L. Pedersen, and A.~W{\k{a}}sowski.
\newblock {Decision problems for interval Markov chains}.
\newblock In {\em Language and Automata Theory and Applications: 5th
  International Conference, LATA 2011, Tarragona, Spain, May 26-31, 2011.
  Proceedings 5}, pages 274--285. Springer, 2011.
\newblock \href {https://doi.org/10.1007/978-3-642-21254-3\_21}
  {\path{doi:10.1007/978-3-642-21254-3\_21}}.

\bibitem{eftang2013port}
J.~L. Eftang and A.~T. Patera.
\newblock Port reduction in parametrized component static condensation:
  approximation and a posteriori error estimation.
\newblock {\em International Journal for Numerical Methods in Engineering},
  96(5):269--302, 2013.
\newblock \href {https://doi.org/10.1002/nme.4543}
  {\path{doi:10.1002/nme.4543}}.

\bibitem{givan2000bounded}
R.~Givan, S.~Leach, and T.~Dean.
\newblock {Bounded-parameter Markov decision processes}.
\newblock {\em Artificial Intelligence}, 122(1-2):71--109, 2000.
\newblock \href {https://doi.org/10.1016/S0004-3702(00)00047-3}
  {\path{doi:10.1016/S0004-3702(00)00047-3}}.

\bibitem{grieves2017digital}
M.~Grieves and J.~Vickers.
\newblock Digital twin: Mitigating unpredictable, undesirable emergent behavior
  in complex systems.
\newblock {\em Transdisciplinary perspectives on complex systems: New findings
  and approaches}, pages 85--113, 2017.
\newblock \href {https://doi.org/10.1007/978-3-319-38756-7\_4}
  {\path{doi:10.1007/978-3-319-38756-7\_4}}.

\bibitem{hahn2011synthesis}
E.~M. Hahn, T.~Han, and L.~Zhang.
\newblock {Synthesis for PCTL in parametric Markov decision processes}.
\newblock In {\em {NASA Formal Methods Symposium}}, pages 146--161. Springer,
  2011.
\newblock \href {https://doi.org/10.1007/978-3-642-20398-5\_12}
  {\path{doi:10.1007/978-3-642-20398-5\_12}}.

\bibitem{holton2003value}
G.~A. Holton.
\newblock {\em {Value-at-Risk: Theory and Practice}}.
\newblock {Academic Press, Elsevier}, 2003.

\bibitem{huynh2013static}
D.~B.~P. Huynh, D.~J. Knezevic, and A.~T. Patera.
\newblock A static condensation reduced basis element method: approximation and
  a posteriori error estimation.
\newblock {\em ESAIM: Mathematical Modelling and Numerical Analysis},
  47(1):213--251, 2013.
\newblock \href {https://doi.org/10.1051/m2an/2012022}
  {\path{doi:10.1051/m2an/2012022}}.

\bibitem{kapteyn2020data}
M.~G. Kapteyn, D.~J. Knezevic, D.~Huynh, M.~Tran, and K.~E. Willcox.
\newblock Data-driven physics-based digital twins via a library of
  component-based reduced-order models.
\newblock {\em International Journal for Numerical Methods in Engineering},
  2020.
\newblock \href {https://doi.org/10.1002/nme.6423}
  {\path{doi:10.1002/nme.6423}}.

\bibitem{kapteyn2020toward}
M.~G. Kapteyn, D.~J. Knezevic, and K.~Willcox.
\newblock Toward predictive digital twins via component-based reduced-order
  models and interpretable machine learning.
\newblock In {\em AIAA Scitech 2020 Forum}, 2020.
\newblock \href {https://doi.org/10.2514/6.2020-0418}
  {\path{doi:10.2514/6.2020-0418}}.

\bibitem{kapteyn2021probabilistic}
M.~G. Kapteyn, J.~V. Pretorius, and K.~E. Willcox.
\newblock A probabilistic graphical model foundation for enabling predictive
  digital twins at scale.
\newblock {\em Nature Computational Science}, 1(5):337--347, 2021.
\newblock \href {https://doi.org/10.1038/s43588-021-00069-0}
  {\path{doi:10.1038/s43588-021-00069-0}}.

\bibitem{kapteyn2022sensing}
M.~G. Kapteyn and K.~E. Willcox.
\newblock {Design of Digital Twin Sensing Strategies via Predictive Modeling
  and Interpretable Machine Learning}.
\newblock {\em Journal of Mechanical Design}, 144(9):091710, 08 2022.
\newblock \href {https://doi.org/10.1115/1.4054907}
  {\path{doi:10.1115/1.4054907}}.

\bibitem{kochenderfer2012next}
M.~J. Kochenderfer, J.~E. Holland, and J.~P. Chryssanthacopoulos.
\newblock {Next-Generation Airborne Collision Avoidance System}.
\newblock {\em Lincoln Laboratory Journal}, 19(1):17--33, 2012.

\bibitem{koller2009probabilistic}
D.~Koller and N.~Friedman.
\newblock {\em {Probabilistic Graphical Models: Principles and Techniques}}.
\newblock MIT Press, Cambridge, MA, 2009.

\bibitem{kwiatkowska2017probabilistic}
M.~Kwiatkowska, G.~Norman, and D.~Parker.
\newblock Probabilistic model checking: Advances and applications.
\newblock In {\em Formal System Verification}, pages 73--121. Springer
  International Publishing, jun 2017.
\newblock URL: \url{https://doi.org/10.1007%2F978-3-319-57685-5_3}, \href
  {https://doi.org/10.1007/978-3-319-57685-5_3}
  {\path{doi:10.1007/978-3-319-57685-5_3}}.

\bibitem{mcclellan2022physics}
A.~McClellan, J.~Lorenzetti, M.~Pavone, and C.~Farhat.
\newblock A physics-based digital twin for model predictive control of
  autonomous unmanned aerial vehicle landing.
\newblock {\em {Philosophical Transactions of the Royal Society A}},
  380(2229):20210204, 2022.
\newblock \href {https://doi.org/10.1098/rsta.2021.0204}
  {\path{doi:10.1098/rsta.2021.0204}}.

\bibitem{morato2023inference}
P.~G. Morato, C.~P. Andriotis, K.~G. Papakonstantinou, and P.~Rigo.
\newblock Inference and dynamic decision-making for deteriorating systems with
  probabilistic dependencies through bayesian networks and deep reinforcement
  learning.
\newblock {\em Reliability Engineering \& System Safety}, 235:109144, 2023.
\newblock \href {https://doi.org/10.1016/j.ress.2023.109144}
  {\path{doi:10.1016/j.ress.2023.109144}}.

\bibitem{murphy2012machine}
K.~P. Murphy.
\newblock {\em {Machine Learning: A Probabilistic Perspective}}.
\newblock MIT Press, Cambridge, MA, 2012.

\bibitem{niederer2021scaling}
S.~A. Niederer, M.~S. Sacks, M.~Girolami, and K.~E. Willcox.
\newblock Scaling digital twins from the artisanal to the industrial.
\newblock {\em Nature Computational Science}, 1(5):313--320, 2021.
\newblock \href {https://doi.org/10.1038/s43588-021-00072-5}
  {\path{doi:10.1038/s43588-021-00072-5}}.

\bibitem{nilim2005robust}
A.~Nilim and L.~El~Ghaoui.
\newblock Robust control of markov decision processes with uncertain transition
  matrices.
\newblock {\em Operations Research}, 53(5):780--798, 2005.
\newblock \href {https://doi.org/10.1287/opre.1050.0216}
  {\path{doi:10.1287/opre.1050.0216}}.

\bibitem{pritsker1997evaluating}
M.~Pritsker.
\newblock {Evaluating Value at Risk Methodologies: Accuracy versus
  Computational Time}.
\newblock {\em Journal of Financial Services Research}, 12(2-3):201--242, 1997.
\newblock \href {https://doi.org/10.1023/A:1007978820465}
  {\path{doi:10.1023/A:1007978820465}}.

\bibitem{puggelli2013polynomial}
A.~Puggelli, W.~Li, A.~L. Sangiovanni-Vincentelli, and S.~A. Seshia.
\newblock {Polynomial-time verification of PCTL properties of MDPs with convex
  uncertainties}.
\newblock In {\em Computer Aided Verification: 25th International Conference,
  CAV 2013, Saint Petersburg, Russia, July 13-19, 2013. Proceedings 25}, pages
  527--542. Springer, 2013.
\newblock \href {https://doi.org/10.1007/978-3-642-39799-8\_35}
  {\path{doi:10.1007/978-3-642-39799-8\_35}}.

\bibitem{puterman2014markov}
M.~L. Puterman.
\newblock {\em {Markov Decision Processes: Discrete Stochastic Dynamic
  Programming}}.
\newblock John Wiley \& Sons, 2014.

\bibitem{rasheed2020digital}
A.~Rasheed, O.~San, and T.~Kvamsdal.
\newblock {Digital twin: Values, challenges and enablers from a modeling
  perspective}.
\newblock {\em IEEE Access}, 8:21980--22012, 2020.
\newblock \href {https://doi.org/10.1109/ACCESS.2020.2970143}
  {\path{doi:10.1109/ACCESS.2020.2970143}}.

\bibitem{reich2015probabilistic}
S.~Reich and C.~Cotter.
\newblock {\em {Probabilistic Forecasting and Bayesian Data Assimilation}}.
\newblock {Cambridge University Press}, 2015.
\newblock \href {https://doi.org/10.1017/CBO9781107706804}
  {\path{doi:10.1017/CBO9781107706804}}.

\bibitem{rockafellar2013superquantiles}
R.~T. Rockafellar and J.~O. Royset.
\newblock {Superquantiles and Their Applications to Risk, Random Variables, and
  Regression}.
\newblock In {\em {Theory Driven by Influential Applications}}, chapter~8,
  pages 151--167. INFORMS, 2013.
\newblock \href {https://doi.org/10.1287/educ.2013.0111}
  {\path{doi:10.1287/educ.2013.0111}}.

\bibitem{rockafellar2000optimization}
R.~T. Rockafellar and S.~Uryasev.
\newblock Optimization of conditional value-at-risk.
\newblock {\em Journal of Risk}, 2(3):21--41, 2000.
\newblock \href {https://doi.org/10.21314/JOR.2000.038}
  {\path{doi:10.21314/JOR.2000.038}}.

\bibitem{rockafellar2002conditional}
R.~T. Rockafellar and S.~Uryasev.
\newblock Conditional value-at-risk for general loss distributions.
\newblock {\em Journal of Banking \& Finance}, 26(7):1443--1471, 2002.
\newblock \href {https://doi.org/10.1016/S0378-4266(02)00271-6}
  {\path{doi:10.1016/S0378-4266(02)00271-6}}.

\bibitem{royset2017risk}
J.~Royset, L.~Bonfiglio, G.~Vernengo, and S.~Brizzolara.
\newblock Risk-adaptive set-based design and applications to shaping a
  hydrofoil.
\newblock {\em Journal of Mechanical Design}, 139(10):101403, 2017.
\newblock \href {https://doi.org/10.1115/1.4037623}
  {\path{doi:10.1115/1.4037623}}.

\bibitem{russell2010artificial}
S.~J. Russell.
\newblock {\em {Artificial Intelligence: A Modern Approach}}.
\newblock Pearson Education, London, UK, 2020.

\bibitem{smetana2016optimal}
K.~Smetana and A.~T. Patera.
\newblock Optimal local approximation spaces for component-based static
  condensation procedures.
\newblock {\em SIAM Journal on Scientific Computing}, 38(5):A3318--A3356, 2016.
\newblock \href {https://doi.org/10.1137/15M1009603}
  {\path{doi:10.1137/15M1009603}}.

\bibitem{sutton2018reinforcement}
R.~S. Sutton and A.~G. Barto.
\newblock {\em {Reinforcement Learning: An Introduction}}.
\newblock MIT press, 2018.

\bibitem{torzoni2024digital}
M.~Torzoni, M.~Tezzele, S.~Mariani, A.~Manzoni, and K.~E. Willcox.
\newblock {A digital twin framework for civil engineering structures}.
\newblock {\em Computer Methods in Applied Mechanics and Engineering},
  418:116584, January 2024.
\newblock \href {https://doi.org/10.1016/j.cma.2023.116584}
  {\path{doi:10.1016/j.cma.2023.116584}}.

\bibitem{tuegel2011reengineering}
E.~J. Tuegel, A.~R. Ingraffea, T.~G. Eason, and S.~M. Spottswood.
\newblock {Reengineering Aircraft Structural Life Prediction Using a Digital
  Twin}.
\newblock {\em International Journal of Aerospace Engineering}, 2011:154798,
  2011.
\newblock \href {https://doi.org/10.1155/2011/154798}
  {\path{doi:10.1155/2011/154798}}.

\bibitem{wiesemann2013robust}
W.~Wiesemann, D.~Kuhn, and B.~Rustem.
\newblock {Robust Markov decision processes}.
\newblock {\em Mathematics of Operations Research}, 38(1):153--183, 2013.
\newblock \href {https://doi.org/10.1287/moor.1120.0566}
  {\path{doi:10.1287/moor.1120.0566}}.

\bibitem{wikle2007bayesian}
C.~K. Wikle and L.~M. Berliner.
\newblock {A Bayesian tutorial for data assimilation}.
\newblock {\em Physica D: Nonlinear Phenomena}, 230(1-2):1--16, 2007.
\newblock \href {https://doi.org/10.1016/j.physd.2006.09.017}
  {\path{doi:10.1016/j.physd.2006.09.017}}.

\bibitem{wolff2012robust}
E.~M. Wolff, U.~Topcu, and R.~M. Murray.
\newblock {Robust control of uncertain Markov decision processes with temporal
  logic specifications}.
\newblock In {\em 2012 IEEE 51st IEEE Conference on Decision and Control
  (CDC)}, pages 3372--3379. IEEE, 2012.
\newblock \href {https://doi.org/10.1109/CDC.2012.6426174}
  {\path{doi:10.1109/CDC.2012.6426174}}.

\bibitem{yang2017algorithms}
H.~Yang and M.~Gunzburger.
\newblock Algorithms and analyses for stochastic optimization for turbofan
  noise reduction using parallel reduced-order modeling.
\newblock {\em Computer Methods in Applied Mechanics and Engineering},
  319:217--239, 2017.
\newblock \href {https://doi.org/10.1016/j.cma.2017.02.030}
  {\path{doi:10.1016/j.cma.2017.02.030}}.

\bibitem{ye2020digital}
Y.~Ye, Q.~Yang, F.~Yang, Y.~Huo, and S.~Meng.
\newblock Digital twin for the structural health management of reusable
  spacecraft: a case study.
\newblock {\em Engineering Fracture Mechanics}, 234:107076, 2020.
\newblock \href {https://doi.org/10.1016/j.engfracmech.2020.107076}
  {\path{doi:10.1016/j.engfracmech.2020.107076}}.

\bibitem{zhang2016decomposition}
W.~Zhang, H.~Rahimian, and G.~Bayraksan.
\newblock {Decomposition Algorithms for Risk-Averse Multistage Stochastic
  Programs with Application to Water Allocation under Uncertainty}.
\newblock {\em INFORMS Journal on Computing}, 28(3):385--404, 2016.
\newblock \href {https://doi.org/10.1287/ijoc.2015.0684}
  {\path{doi:10.1287/ijoc.2015.0684}}.

\end{thebibliography}
